\documentclass[11pt,a4paper,oneside]{amsart}

\usepackage[a4paper,inner=2.5cm,outer=2.5cm,top=3cm,bottom=3cm,pdftex]{geometry}

\usepackage{color}
\usepackage{amsmath,amssymb,amsfonts}
\usepackage{graphicx}

\usepackage[utf8]{inputenc}
\usepackage[colorlinks=true,linkcolor=black,citecolor=black]{hyperref}

\newtheorem{thm}{Theorem}[section]

\newtheorem{lem}[thm]{Lemma}

\newtheorem{defn}[thm]{Definition}

\theoremstyle{plain}
\newtheorem{rem}[thm]{Remark}

\numberwithin{equation}{section}
\newcommand{\RR}{\mathbb{R}}

\newcommand{\TT}{\mathbb{T}}

\newcommand{\ZZ}{\mathbb{Z}}
\newcommand{\NN}{\mathbb{N}}
\newcommand{\PP}{\mathbb{P}}

\newcommand{\DD}{\Delta}

\newcommand{\Bcal}{\mathcal{B}}
\newcommand{\Ccal}{\mathcal{C}}
\newcommand{\Fcal}{\mathcal{F}}
\newcommand{\Hcal}{\mathcal{H}}
\newcommand{\Kcal}{\mathcal{K}}
\newcommand{\Lcal}{\mathcal{L}}
\newcommand{\Ocal}{\mathcal{O}}

\newcommand{\Scal}{\mathcal{S}}
\newcommand{\Ucal}{\mathcal{U}}
\newcommand{\Vcal}{\mathcal{V}}
\newcommand{\Wcal}{\mathcal{W}}
\newcommand{\Zcal}{\mathcal{Z}}

\newcommand{\ee}{\varepsilon}

\newcommand{\dd}{\partial}

\newcommand{\divv}{\mbox{div\,}}
\newcommand{\rot}{\mbox{rot\,}}

\newcommand{\abs}[1]{\left\vert#1\right\vert}
\newcommand{\set}[1]{\left\{#1\right\}}
\newcommand{\psca}[1]{\left\langle#1\right\rangle}
\newcommand{\pint}[1]{\left[#1\right]}
\newcommand{\pare}[1]{\left(#1\right)}
\newcommand{\norm}[1]{\left\Vert#1\right\Vert}

\newcommand{\ds}{\displaystyle}

\begin{document}

\title{Approximate controllability of second grade fluids}%
\author{Van-Sang Ngo}
\address{Laboratoire de Math\'ematiques Rapha\"el Salem, UMR 6085 CNRS, Universit\'e de Rouen, 76801 Saint-Etienne du Rouvray Cedex, France}
\email{van-sang.ngo@univ-rouen.fr}
\author{Genevi\`eve Raugel}
\address{CNRS, Laboratoire de Math\'ematiques d'Orsay, Université de Paris-Sud, Orsay Cedex, F-91405, France}
\email{genevieve.raugel@math.u-psud.fr}

%\thanks{This work was partially done during the visit of the author at the Vietnam Institute for Advanced Study in Mathematics (VIASM). The author is thankful to the VIASM for the financial support and for the very kind hospitality of the institute and of all the staff.}

\begin{abstract}
    This paper deals with the controllability of the second grade fluids, a class of non-Newtonian of differentiel type, on a two-dimensional torus. Using the method of Agrachev-Sarychev \cite{AS2005}, \cite{AS2006} and of Sirikyan \cite{S2006}, we prove that the system of second grade fluids is approximately controllable by a finite-dimensional control force.
\end{abstract}

\maketitle

%----- %----- %----- %----- %----- %

\bigskip
\section{Introduction} \label{se:Intro}

The goal of this paper is to study the approximate controllability of the system of fluids of second grade, using low-mode (finite-dimensional) control forces. More precisely, we consider the following system
\begin{equation}
	\label{sys:FG2eta}
	\left\{
	\begin{aligned}
		&\dd_t\pare{u-\alpha\DD u} - \nu\DD u + \rot\!\pare{u-\alpha\DD u}\times u + \nabla p = f + \eta\\
		&\divv u = 0\\
		&u(0) = u_0,
	\end{aligned}
	\right.
\end{equation}
on the domain $\TT^2$, which is the two-dimensional torus $]0,2\pi q_1[ \times ]0,2\pi q_2[$, with $q_1 > 0$ and $q_2 > 0$. Here $u = (u_1(t,x),u_2(t,x))$ and $p = p(t,x)$ are unknown and represent the velocity vector field and the pressure function; $f = f(t,x)$ is the extenal force field; and the control force $\eta = \eta(t,x)$ is supposed to belong to a finite-dimensional space which will be made more precise later.

Fluids of second grade belong to a particular class of non-Newtonian Rivlin-Ericksen fluids of differential type \cite{RE}, which usually arise in petroleum industry, in polymer technology or in liquid crystal suspension problems. For these fluids, the Cauchy stress tensor $\sigma$ is not linearly proportional to the local strain rate but given by
\begin{equation}
	\label{eq:StressTensor} \sigma = -pI + 2\nu A_1 + \alpha_1 A_2 + \alpha_2 A_1^2,
\end{equation}
where $\nu$ stands for the kinematic viscosity, $p$ is the pressure and $A_1$, $A_2$ represent the first two Rivlin-Ericksen tensors, which are
\begin{equation*}
	A_1(u) = \frac{1}{2} \pare{\nabla u + \nabla u^T},
\end{equation*}
corresponding to the local strain tensor and
\begin{equation*}
	A_2(u) = \frac{DA_1}{Dt} + \pare{\nabla u}^T A_1 + A_1 \pare{\nabla u},
\end{equation*}
where 
\begin{equation*}
	\frac{D}{Dt} = \dd_t + u\cdot\nabla
\end{equation*}
is the material derivative. In \cite{DF}, Dunn and Fosdick used the compatibility of \eqref{eq:StressTensor} with thermodynamics to prove that
\begin{equation*}
	\alpha_1 + \alpha_2 = 0; \quad \alpha_1 \geq 0.
\end{equation*}
Setting $\alpha = \alpha_1$ and writing the equation
\begin{equation*}
	\frac{Du}{Dt} = \dd_t u + u\cdot\nabla u = \divv \sigma + f,
\end{equation*}
one obtain the equations of second grade fluids of the following form
\begin{equation}
	\label{sys:FG2}
	\left\{
	\begin{aligned}
		&\dd_t\pare{u-\alpha\DD u} - \nu\DD u + \rot\!\pare{u-\alpha\DD u}\times u + \nabla p = f && \mbox{in } \RR_+\times \TT^2\\
		&\divv u = 0 && \mbox{in } \RR_+\times \TT^2\\
		&u(0) = u_0 && \mbox{in } \TT^2.
	\end{aligned}
	\right.
\end{equation}

The local existence in time and uniqueness of a strong solution to \eqref{sys:FG2} have been proven by Cioranescu and Ouazar in \cite{CO} in the case of two-dimensional or three-dimensional domains with non-slip boundary conditions. Moreover, the solution is global in time in the two-dimensional case. Second grade fluids in these domains were also studied by Moise, Rosa and Wang in \cite{MRW}, where the authors proved the existence of a compact global attractor in the two-dimensional case. The existence, the uniqueness of a strong solution and the dynamics of second grade fluids in the torus $\TT^2$ was studied in \cite{PRR} by Paicu, Raugel and Rekalo, and in \cite{PR} by the first two authors, using the Lagrangian approach. For further results concerning the system \eqref{sys:FG2}, we refer the readers to \cite{B}, \cite{BL}, \cite{CG}, \cite{GaR}, \cite{GaS}, \cite{GaGS1}, \cite{GaGS2}, \cite{GSa}, \cite{GSc}, \cite{LR}, \ldots

In this paper, in order to study the approximate controllability of the second grade fluid system \eqref{sys:FG2eta} by a low-mode control $\eta$, we use the method introduced by Agrachev and Sarychev in \cite{AS2005} and \cite{AS2006} for the Navier-Stokes and Euler systems in the two-dimensional torus $\TT^2$. This method was extended later for the three-dimensional Navier-Stokes system by Shirikyan in \cite{S2006} and \cite{S2007} and for the three-dimensional Euler system by Nersisyan in \cite{N}. The main idea consists in proving that, if the (finite-dimensional) space of controls $E$ contains sufficiently many Fourier modes then, for any $T > 0$, the system \eqref{sys:FG2eta} is approximately controllable in time $T$ by an $E$-valued control $\eta$.

\smallskip

Before stating the main results and the main ideas of this paper, we will introduce the needed notations and function spaces. Let $H^m(\TT^2)^2$ be the classical Sobolev space of two-dimensional vector fields, whose components belong to $H^m(\TT^2)$. For $m=0$, we simply have $H^0(\TT^2)^2 = L^2(\TT^2)^2$. As in \cite{PRR}, for any $m \in \NN$, we denote $V^m(\TT^2)^2$ the closure of the space
\begin{equation*}
	\set{u \in C^\infty(\TT^2)^2 \;\vert\; u \mbox{ is periodic }, \divv u = 0, \int_{\TT^2} u dx = 0}
\end{equation*}
in $H^m(\TT^2)^2$. Then $V^m(\TT^2)^2$ is a Banach space, endowed with the classical norm of $H^m(\TT^2)^2$. For any $\theta > 0$, we define $V^\theta(\TT^2)^2$ using the method of interpolation between $V^m(\TT^2)^2$ spaces. Finally, we also use $H^m_{per}(\TT^2)^2$ to denote the space of vector fields $u \in H^m(\TT^2)^2$, which are periodic and whose mean value is zero. 

In what follows, we recall the definition of a strong solution of the system \eqref{sys:FG2}.
\begin{defn}
	\label{de:StrongSol}
	Let $T > 0$. For any $f \in L^\infty\pare{0,T,H^1_{per}(\TT^2)^2}$ and $u_0 \in V^3(\TT^2)^2$, the vector field $u(t,x)$ is said to be a strong solution of the system \eqref{sys:FG2}, with data $(f,u_0)$, on the time interval $[0,T]$ if $u~\in~C(0,T,V^3(\TT^2)^2)$, $\dd_t u \in L^\infty(0,T,V^2(\TT^2)^2)$, $u(0) = u_0$, and for any $t \in ]0,T]$, for any $\phi \in V^0(\TT^2)^2$, the following equation holds
	\begin{equation}
		\label{eq:StrongSol} \psca{\dd_t\pare{u(t)-\alpha\DD u(t)} - \nu\DD u(t) + \rot\!\pare{u(t)-\alpha\DD u(t)} \times u(t), \phi} = \psca{f(t),\phi}.
	\end{equation}
\end{defn}
\noindent In \cite{PRR}, the authors prove that
\begin{thm}
	Let $\alpha > 0$ and $T > 0$. 
	\begin{enumerate}
		\item For any $f \in L^\infty(0,T,H^1_{per}(\TT^2)^)$ and any $u_0 \in V^3(\TT^2)^2$, there exists a unique strong solution $$u \in C(0,T,V^3(\TT^2)^2) \cap W^{1,\infty}(0,T,V^2(\TT^2)^2)$$ of the system \eqref{sys:FG2}. Moreover, for any $t\in[0,T]$, the map $$V^3(\TT^2)^2 \ni u_0 \mapsto u(t) \in V^3(\TT^2)^2$$ is continuous.
		\item Let $m \geq 1$. Assume that $f\in L^\infty(0,T,H^{m+1}_{per}(\TT^2)^2)$ and $u_0 \in V^{m+3}(\TT^2)^2$. Then, the solution $u$ of the system \eqref{sys:FG2} belongs to $C(0,T,V^{m+3}(\TT^2)^2)$.
	\end{enumerate}
\end{thm}

For the system \eqref{sys:FG2eta}, we want to define the approximate controllability using low-mode controls. We will adapt the definition of approximate controllability given in \cite{S2006} to the case of fluids of second grade.
\begin{defn}
	\label{de:ApproxCont}
	Let $\theta > 0$, $T > 0$ and let $E$ be a finite-dimensional subspace of $V^3(\TT^2)^2$. The second grade fluid system \eqref{sys:FG2eta} is $\theta$-approximately controllable ($\theta$-AC) in time $T$ by $E$-valued controls if, for any $\ee > 0$, for any $u_0, u_T \in V^3(\TT^2)$, there exist a control $\eta \in L^\infty(0,T,E)^2$ and a strong solution $u \in C(0,T,V^3(\TT^2)^2)$ of the system \eqref{sys:FG2eta} such that
	\begin{equation*}
		\norm{u(T) - u_T}_{V^\theta(\TT^2)^2} \leq \ee.
	\end{equation*}
\end{defn}

For any $m \in \ZZ^2 \setminus \set{0}$, let $$c_m(x) = m^{q,\perp}\cos\psca{m,x}_q \qquad \mbox{ and } \qquad s_m(x) = m^{q,\perp}\sin\psca{m,x}_q,$$ where $m^{q,\perp}$ will be defined in Section \ref{se:Saturation}. Then, it is classical that $c_m$, $s_m$, with $m \in \ZZ^2 \setminus \set{0}$, are eigenvectors of the Stokes operator $-\mathbb{P}\DD$, where $\mathbb{P}$ is the Leray projection, and that the family $\set{c_m,s_m \;\vert\; m \in \ZZ^2 \setminus \set{0}}$ forms an orthonormal basis of $V^k(\TT^2_q)^2$, $k\in \NN$. For any $N\in\NN^*$, we set 
\begin{equation}
	\label{eq:HN} \Hcal^N_q = Span\set{c_m,s_m \;\vert\; m\in \ZZ\setminus\set{0}, \abs{m} \leq N}.
\end{equation}
The main result of this paper is the following theorem.
%\begin{thm}
%	\label{th:FG2approx} Let $0 \leq \theta < 3$, $T > 0$ and $f\in L^\infty(0,T,H^2_{per}(\TT^2)^2)$. Then the system \eqref{sys:FG2eta} is $\theta$-AC in time $T$ by $\Hcal^3_q$-valued controls. Moreover, if $u_0, u_T \in V^4(\TT^2)$ then the system \eqref{sys:FG2eta} is $3$-AC in time $T$ by $\Hcal^3_q$-valued controls.
%\end{thm}
\begin{thm}
	\label{th:FG2approx} Let $T > 0$, $f\in L^\infty(0,T,H^2_{per}(\TT^2)^2)$ and $u_0, u_T \in V^4(\TT^2)$. Then the system \eqref{sys:FG2eta} is $3$-AC in time $T$ by $\Hcal^3_q$-valued controls.
\end{thm}

We note that, unlike the case of Navier-Stokes equations, the system of second grade fluids  is an exemple of asymptotically smooth system, which only possesses a smoothing effect in infinite time. The systems \eqref{sys:FG2eta} or \eqref{sys:FG2} also differ from the $\alpha$-type models, such as the so-called $\alpha$-Navier-Stokes system (see \cite{FHT01}, \cite{FHT02} and the references therein). Indeed, the $\alpha$-Navier-Stokes system contains the very regularizing term $-\nu\DD (u-\alpha\DD u)$ instead of $-\nu\DD u$, and thus is a semilinear problem, which is easier to solve. It is different in the case of second grade fluids where the dissipation is much weaker. This weak smoothing effect explains why in our result, we can not obtain an approximate control in the same norm as the initial data but only a control in the weaker norm. We also remark the similar phenomenon in \cite{PRR}, where the Navier-Stokes system is proved to be the limit of the second grade fluid system in $V^\theta(\TT^2)^2$ for data in $V^3(\TT^2)^2$ (respectively in $V^3(\TT^2)^2$ for data in $V^4(\TT^2)^2$).

Another problem when we want to apply the method of Shirikyan \cite{S2006} lies in the complexity of the nonlinear term and the appearance of the term $\dd_t \pare{-\alpha\DD u}$. To avoid this difficulty, let $$\Ucal = u-\alpha\DD u \quad\mbox{ and }\quad \Ucal_0 = u_0-\alpha\DD u_0$$ and let us rewrite the system \eqref{sys:FG2eta} in the following form
\begin{equation}
	\label{sys:FG2etaCal}
	\left\{
	\begin{aligned}
		&\dd_t \Ucal + \Lcal\Ucal + \Bcal(\Ucal,\Ucal) = \PP f + \eta\\
		&\divv \Ucal = 0\\
		&\Ucal(0) = \Ucal_0,
	\end{aligned}
	\right.
\end{equation}
where
\begin{equation}
	\label{eq:LcalBcal}
	\left\{ 
	\begin{aligned}
		&\Lcal\Ucal = - \nu\PP \DD (I-\alpha\DD)^{-1} \Ucal\\
		&\Bcal(\Ucal_1,\Ucal_2) = \PP \pare{\rot \Ucal_1\times \pare{(I-\alpha\DD)^{-1} \Ucal_2}}.
	\end{aligned}
	\right.
\end{equation}
Along with the system \eqref{sys:FG2etaCal}, we consider the following controlled system
\begin{equation}
	\label{sys:FG2etazetaCal}
	\left\{
	\begin{aligned}
		&\dd_t \Ucal + \Lcal(\Ucal+\zeta) + \Bcal(\Ucal+\zeta,\Ucal+\zeta) = \PP f + \eta\\
		&\divv \Ucal = 0\\
		&\Ucal(0) = \Ucal_0.
	\end{aligned}
	\right.
\end{equation}
As in \cite{AS2005} or \cite{S2006}, we give the following definition
\begin{defn}
	\label{de:Convexification} For any finite-dimensional subspace $E$ of $V^3(\TT^2)^2$, we define $\Fcal(E)$ as the largest vector subspace of $V^3(\TT^2)^2$ such that, $\Fcal(E) \supset E$, and if $\overline{\eta} \in \Fcal(E) \setminus E$ then, there exist $$k\in\NN^*,\; \alpha_1,\ldots, \alpha_k > 0,\; \eta, \rho^1, \ldots, \rho^k \in E$$ satisfying
	\begin{equation*}
		\overline{\eta} = \eta - \sum_{j=1}^k \alpha_j \Bcal(\rho^j,\rho^j). 
	\end{equation*}
\end{defn}
The main idea to prove Theorem \ref{th:FG2approx} is to extend the space of control $E$ to the larger space $\Fcal(E)$. To this end, we will prove the following theorems.
\begin{thm}
	\label{th:FG2extension} Let $T > 0$ and $E$ a finite-dimensional subspace of $V^3(\TT^2)^2$. Then, the system \eqref{sys:FG2etaCal} is ($\theta$-) approximately controllable in time $T$ by an $E$-valued control $\eta$ if and only if so is the system \eqref{sys:FG2etazetaCal} with $E$-valued controls $\eta$ and $\zeta$.
\end{thm}
\begin{thm}
	\label{th:FG2convex} Let $T > 0$ and $E$ a finite-dimensional subspace of $V^3(\TT^2)^2$. Then, the system \eqref{sys:FG2etazetaCal} is approximately controllable in time $T$ by $E$-valued controls $\eta$ and $\zeta$ if and only if so is the system \eqref{sys:FG2etazetaCal} with $\Fcal(E)$-valued controls $\eta$.
\end{thm}
As a consequence of these theorems, the approximate controllability by $E$-valued controls is equivalent to the approximate controllability by $\Fcal(E)$-valued controls. We can define a sequence of subspace $$E = E_0 \subset E_1 \subset \ldots \subset E_n \subset \ldots$$ such that, for any $n\in\NN$ we have $E_{n+1} = \Fcal(E_n)$. The only problem is that $\Fcal(E)$ may be not larger than $E$. However, if we can choose $E$ in such a way that the space $$E_\infty = \bigcup_{n=0}^\infty E_n$$ is dense in $V^1(\TT^2)^2$, then the approximate controllability $E$-valued controls will follows. Indeed, for $T > 0$, $\ee > 0$, $u_0, u_T \in V^4(\TT^2)^2$, we set $$\Ucal_0 = u_0 - \alpha\DD u_0 \quad\mbox{and}\quad \Ucal_T = u_T - \alpha\DD u_T.$$ Then, we can always exactly control the system \eqref{sys:FG2etaCal} by the control 
$$\eta = \dd_t\overline{\Ucal} + \Lcal\overline{\Ucal} + \Bcal(\overline{\Ucal}) - \PP f \in L^{\infty}(0,T,V^1(\TT^2)^2),$$ where 
$$\overline{\Ucal}(t) = \frac{1}T \pare{I-\alpha\DD} \pare{(T-t)u_0 + tu_T}.$$ Thus, if $E_\infty$ is dense in $V^1(\TT^2)^2$, then there exists $n \in \NN$ large enough such that the system \eqref{sys:FG2etaCal} is approximately controllable by $E_n$-valued controls. Using Theorems \ref{th:FG2extension} and \ref{th:FG2convex}, we can prove by induction that \eqref{sys:FG2etaCal} is approximately controllable by $E$-valued controls. In order to prove Theorem \ref{th:FG2approx}, we only need to prove the following result
\begin{thm}
	\label{th:FG2Saturation} If $E = \Hcal^3_q$ then $E_\infty \supset \Hcal^N_q$, for any $N \in \NN$, $N \geq 3$.
\end{thm}
The paper will be organized as follows. In Section \ref{se:Pre}, we study a pertubation of the system \eqref{sys:FG2etaCal}, which is necessary to prove our main theorem. Theorem \ref{th:FG2extension} will be proved in Section \ref{se:Extension}. Section \ref{se:Convexification} is devoted to the demonstration of Theorem \ref{th:FG2convex}. In Section \ref{se:Saturation}, we put in evidence the saturation property given in Theorem \ref{th:FG2Saturation}. Finally, in the last section, we wil prove the main theorem \ref{th:FG2approx}.

%----- %----- %----- %----- %----- %

\bigskip
\section{Preliminary results on the system of fluids of second grade} \label{se:Pre}

In this section, we consider the following perturbed system of fluids of second grade
\begin{equation}
    \label{sys:FG2perCal}
    \left\{
    \begin{aligned}
        &\partial_t \Wcal + \Lcal\Wcal + \Bcal(\Wcal) + \Bcal(\Wcal,\Vcal) + \Bcal(\Vcal,\Wcal) = \PP f &&\mbox{in } \TT^2 \times [0,T]\\
        &\divv \Wcal = 0 &&\mbox{in } \TT^2 \times [0,T]\\
        &\Wcal(0) = \Wcal_0 &&\mbox{in } \TT^2,
    \end{aligned}
    \right.
\end{equation}
where $\Vcal \in L^{\infty}(0,T,V^2(\TT^2)^2)$, $f\in L^{\infty}(0,T,H^1_{per}(\TT^2)^2)$. We want to prove that
\begin{thm}
	\label{th:FG2pertCal} Let $T > 0$ fixed.
	
	\noindent 1. For any $\Vcal \in L^\infty(0,T,V^2(\TT^2)^2)$, $f \in L^\infty(0,T,H^1_{per}(\TT^2)^2)$ and $\Wcal_0 \in V^1(\TT^2)^2$, the system \eqref{sys:FG2perCal} has a unique solution $$\Wcal \in L^{\infty}(0,T,V^1(\TT^2)^2) \cap C(0,T,V^1(\TT^2)^2).$$ Moreover, if $\Vcal \in L^\infty(0,T,V^3(\TT^2)^2)$, $f \in L^\infty(0,T,H^2_{per}(\TT^2)^2)$ and $\Wcal_0 \in V^2(\TT^2)^2$ then 
	$$\Wcal \in L^{\infty}(0,T,V^2(\TT^2)^2) \cap C(0,T,V^2(\TT^2)^2).$$
		
	\noindent 2. Let $\Wcal$ and $\widehat{\Wcal}$ be two solutions of the system \eqref{sys:FG2perCal}, corresponding to data $(\Vcal,f,\Wcal_0)$ and $(\widehat{\Vcal},\widehat{f},\widehat{\Wcal}_0)$ respectively. Then, if $$\widehat{\Wcal} \in L^{\infty}(0,T,V^2(\TT^2)^2) \cap C(0,T,V^2(\TT^2)^2),$$ then, there exists a constant $C>0$ such that, for any $t\in [0,T]$, we have
	\begin{multline}
		\label{eq:LipschitzCalbis} \norm{\Wcal(t)-\widehat{\Wcal}(t)}_{V^1(\TT^2)^2}\\ 
		\leq C \pare{\norm{\Vcal-\widehat{\Vcal}}_{L^2(0,T,V^2)^2} + \norm{\PP f-\PP \widehat{f}}_{L^2(0,T,V^1)^2} + \norm{\Wcal_0-\widehat{\Wcal}_0}_{V^1(\TT^2)^2}}.
	\end{multline}
\end{thm}

We remark that if we set 
\begin{equation*}
	\left\{ 
	\begin{aligned}
		v &= (I-\alpha\DD)^{-1} \Vcal\\
		w_0 &= (I-\alpha\DD)^{-1} \Wcal_0\\
		w &= (I-\alpha\DD)^{-1} \Wcal
	\end{aligned}
	\right.
\end{equation*}
then, $w$ is solution of the following system
\begin{equation}
    \label{sys:FG2perP}
    \left\{
    \begin{aligned}
        &\partial_t(w-\alpha\DD w) - \nu\DD w + \PP B(w) + \PP B(w,v) + \PP B(v,w) = \PP f\\
        &\divv w = 0\\
        &w(0) = w_0,
    \end{aligned}
    \right.
\end{equation}
where
$$B(u_1,u_2) = \rot(u_1-\alpha\DD u_1) \times u_2 \qquad\mbox{and}\qquad B(u) = B(u,u).$$
Theorem \ref{th:FG2pertCal} is in fact equivalent to the following theorem for the system \eqref{sys:FG2perP}
\begin{thm}
	\label{th:FG2pert2D} Let $T > 0$ fixed.
	
	\noindent 1. For any $v \in L^\infty(0,T,V^4(\TT^2)^2)$, $f \in L^\infty(0,T,H^1_{per}(\TT^2)^2)$ and $w_0 \in V^3(\TT^2)^2$, the system \eqref{sys:FG2perP} has a unique solution $$w \in L^{\infty}(0,T,V^3(\TT^2)^2) \cap C(0,T,V^3(\TT^2)^2).$$ Moreover, if $v \in L^\infty(0,T,V^3(\TT^2)^2)$, $f \in L^\infty(0,T,H^2_{per}(\TT^2)^2)$ and $w_0 \in V^2(\TT^2)^2$ then 
	$$w \in L^{\infty}(0,T,V^4(\TT^2)^2) \cap C(0,T,V^4(\TT^2)^2).$$
		
	\noindent 2. Let $w$ and $\widehat{w}$ be two solutions of the system \eqref{sys:FG2perP}, corresponding to data $(v,f,w_0)$ and $(\widehat{v},\widehat{f},\widehat{w}_0)$ respectively. Then, if $$\widehat{w} \in L^{\infty}(0,T,V^4(\TT^2)^2) \cap C(0,T,V^4(\TT^2)^2),$$ then, there exists a constant $C>0$ such that, for any $t\in [0,T]$, we have
	\begin{multline}
		\label{eq:Lipschitz2Dbis} \norm{w(t)-\widehat{w}(t)}_{V^3(\TT^2)^2}\\ 
		\leq C \pare{\norm{v-\widehat{v}}_{L^2(0,T,V^4)^2} + \norm{\PP f-\PP \widehat{f}}_{L^2(0,T,V^1)^2} + \norm{w_0-\widehat{w}_0}_{V^3(\TT^2)^2}}.
	\end{multline}
\end{thm}

\bigskip

We remark that the proof of the existence of solutions of the systems \eqref{sys:FG2perCal} and \eqref{sys:FG2perP} follows the lines of the proof of [\cite{PRR}, Theorems 2.1 and 2.4]. In what follows, we give the needed \emph{a priori} estimates to prove \eqref{eq:LipschitzCalbis} and \eqref{eq:Lipschitz2Dbis}. We will set 
\begin{equation}
	\label{eq:Notations}
	\begin{aligned}
		&\Wcal = \rot(w-\alpha\DD w), \qquad \qquad &&\widehat{\Wcal} = \rot(\widehat{w}-\alpha\DD\widehat{w})\\
		&\Vcal = \rot(v-\alpha\DD v), &&\widehat{\Vcal} = \rot(\widehat{v}-\alpha\DD\widehat{v})\\
		&\Ocal = \Wcal - \widehat{\Wcal}. && %W = w - \widehat{w}.
	\end{aligned}
\end{equation}

%----- %----- %----- %----- %----- %

\bigskip

\subsection{Propagation of the $V^3$-norm}

In this paragraph, we give a priori estimates of a solution of the system \eqref{sys:FG2perP} in $V^3(\TT^2)$-norm. Applying the $\rot$ operator to the first equation of \eqref{sys:FG2perP}, we obtain
\begin{equation}
	\label{eq:propag2D} \dd_t\Wcal + \frac{\nu}{\alpha}\Wcal + \PP \pare{w\cdot\nabla\Wcal} + \PP \pare{v\cdot\nabla\Wcal} + \PP \pare{w\cdot\nabla \Vcal} = \rot \PP f + \frac{\nu}{\alpha} \rot w.
\end{equation}
Since $v$ and $w$ are divergence-free vector fields on $\TT^2$, integrations by parts show that
$$\psca{w\cdot\nabla\Wcal,\Wcal}_{L^2(\TT^2)^2} = \psca{v\cdot\nabla\Wcal,\Wcal}_{L^2(\TT^2)^2} = 0.$$ As a consequence, taking the $L^2(\TT^2)^2$ inner product of \eqref{eq:propag2D} with $\Wcal$, we get
\begin{align}
	\label{eq:propag2DL2} &\frac{1}2\frac{d}{dt} \norm{\Wcal}_{L^2(\TT^2)^2}^2 + \frac{\nu}{\alpha} \norm{\Wcal}_{L^2(\TT^2)^2}^2\\ 
	&\qquad \qquad \qquad \leq \abs{\psca{w\cdot\nabla \Vcal,\Wcal}_{L^2(\TT^2)^2}} + \abs{\psca{\rot f,\Wcal}_{L^2(\TT^2)^2}} + \frac{\nu}{\alpha} \abs{\psca{\rot w,\Wcal}_{L^2(\TT^2)^2}}\notag\\
	&\qquad \qquad \qquad \leq \norm{\rot f}_{L^2(\TT^2)^2} \norm{\Wcal}_{L^2(\TT^2)^2} + \pare{\norm{\nabla \Vcal}_{L^2(\TT^2)^2} + \frac{\nu}{\alpha}} \norm{\Wcal}_{L^2(\TT^2)^2}^2, \notag
\end{align}
which implies that
\begin{equation}
	\label{eq:propag2DL2bis} \frac{d}{dt} \norm{\Wcal}_{L^2(\TT^2)^2}^2 \leq \norm{f}_{V^1(\TT^2)^2}^2 + 2\pare{1 + \norm{v}_{V^4(\TT^2)^2}} \norm{\Wcal}_{L^2(\TT^2)^2}^2.
\end{equation}
Finally, the Gronwall lemma gives, for any $0 \leq t \leq T$,
\begin{multline}
	\label{eq:propag2DL2gron} \norm{\Wcal(t)}_{L^2(\TT^2)^2}^2 \leq \\ \pare{\norm{\Wcal(0)}_{L^2(\TT^2)^2}^2 + \norm{f}_{L^\infty(0,T,V^1(\TT^2)^2)}^2} \exp\set{2t\pare{1 + \norm{v}_{L^\infty(0,T,V^4(\TT^2)^2)}}}.
\end{multline}

%----- %----- %----- %----- %----- %

\bigskip

\subsection{A priori estimates of the difference of two solutions in $V^3$-norm}

Using the same notations as in \cite{CO}, \cite{PRR} or \cite{PR}, we identify a 2D vector $(u_1,u_2)$ with the 3D vector $(u_1,u_2,0)$ and a scalar $\lambda$ with the 3D vector $(0,0,\lambda)$. We also recall the identity $$\rot(\rot(a)\times b) = b\cdot\nabla \rot(a),$$ where $a$ and $b$ are 2D vector fields and $\divv v = 0$. We deduce from \eqref{sys:FG2perP} that $\Wcal$ and $\widehat{\Wcal}$ are solutions of the following equation, with data $(v,f,w_0)$ and $(\widehat{v},\widehat{f},\widehat{w}_0)$ respectively.
\begin{equation}
	\label{eq:rotFG2per2D} \dd_t\Wcal + \frac{\nu}{\alpha}\Wcal + \PP\pare{w\cdot\nabla\Wcal} + \PP\pare{v\cdot\nabla\Wcal} + \PP\pare{w\cdot\nabla \Vcal} = \rot \PP f + \frac{\nu}{\alpha} \rot w.
\end{equation}
The calculation of the difference between the equations corresponding to $\Wcal$ and $\widehat{\Wcal}$ shows that $\Ocal = \Wcal - \widehat{\Wcal}$ satisfies the following equation
\begin{multline}
	\label{eq:Ocal2D} \dd_t\Ocal + \frac{\nu}{\alpha}\Ocal + \PP\pare{w\cdot\nabla\Ocal} + \PP\pare{(w-\widehat{w})\cdot\nabla\widehat{\Wcal}} + \PP\pare{v\cdot\nabla\Ocal} + \PP\pare{(v-\widehat{v})\cdot\nabla\widehat{\Wcal}} \\
	+ \PP\pare{w\cdot\nabla \pare{\Vcal - \widehat{\Vcal}}} + \PP\pare{(w-\widehat{w})\cdot\nabla \widehat{\Vcal}} \;=\; \rot(\PP f - \PP\widehat{f}) + \frac{\nu}{\alpha} \rot(w-\widehat{w}). %\notag
\end{multline}
Next, we will take the $L^2$ inner product of \eqref{eq:Ocal2D} with $\Ocal$.

%----- %----- %----- %----- %----- %

Using the divergence-free property of $v$ and $w$, we have 
\begin{equation}
	\label{eq:Ocal2D01} \psca{v\cdot\nabla\Ocal,\Ocal}_{L^2(\TT^2)^2} = \psca{w\cdot\nabla\Ocal,\Ocal}_{L^2(\TT^2)^2} = 0.
\end{equation}
Now, using H\"older's and Cauchy-Schwarz's inequalities, we get
\begin{align}
	\label{eq:Ocal2D02} \abs{\psca{(w-\widehat{w})\cdot\nabla\widehat{\Wcal},\Ocal}_{L^2(\TT^2)^2}} &\leq \norm{\nabla\widehat{\Wcal}}_{L^2(\TT^2)^2} \norm{w-\widehat{w}}_{L^\infty(\TT^2)^2} \norm{\Ocal}_{L^2(\TT^2)^2}\\
	&\leq C \norm{\nabla\widehat{\Wcal}}_{L^2(\TT^2)^2} \norm{w-\widehat{w}}_{V^3(\TT^2)^2} \norm{\Ocal}_{L^2(\TT^2)^2} \notag\\
	&\leq C \norm{\widehat{w}}_{V^4(\TT^2)^2} \norm{\Ocal}_{L^2(\TT^2)^2}^2. \notag
\end{align}
The same calculations give
\begin{align}
	\label{eq:Ocal2D06} \abs{\psca{(v-\widehat{v})\cdot\nabla\widehat{\Wcal},\Ocal}_{L^2(\TT^2)^2}} &\leq C \norm{\nabla\widehat{\Wcal}}_{L^2(\TT^2)^2} \norm{v-\widehat{v}}_{L^\infty(\TT^2)^2} \norm{\Ocal}_{L^2(\TT^2)^2} \\
	&\leq C \norm{\widehat{w}}_{V^4(\TT^2)^2} \pare{\norm{v-\widehat{v}}_{V^3(\TT^2)^2}^2 + \norm{\Ocal}_{L^2(\TT^2)^2}^2}, \notag
\end{align}
\begin{align}
	\label{eq:Ocal2D09} &\abs{\psca{w\cdot\nabla \pare{\Vcal - \widehat{\Vcal}},\Ocal}_{L^2(\TT^2)^2}}\\ 
	&\leq \norm{w}_{L^\infty(\TT^2)^2} \norm{\nabla\pare{\Vcal - \widehat{\Vcal}}}_{L^2(\TT^2)^2} \norm{\Ocal}_{L^2(\TT^2)^2}\notag \\
	&\leq C \norm{w}_{L^\infty(\TT^2)^2} \norm{v-\widehat{v}}_{V^4(\TT^2)^2} \norm{\Ocal}_{L^2(\TT^2)^2}\notag\\
	&\leq C \pare{\norm{w-\widehat{w}}_{L^\infty(\TT^2)^2} + \norm{\widehat{w}}_{L^\infty(\TT^2)^2}} \norm{v-\widehat{v}}_{V^4(\TT^2)^2} \norm{\Ocal}_{L^2(\TT^2)^2} \notag\\
	&\leq C \pare{\norm{\Ocal}_{L^2(\TT^2)^2} + \norm{\widehat{w}}_{L^\infty(\TT^2)^2}} \norm{v-\widehat{v}}_{V^4(\TT^2)^2} \norm{\Ocal}_{L^2(\TT^2)^2} \notag\\
	&\leq C \norm{\widehat{w}}_{L^\infty(\TT^2)^2} \norm{v-\widehat{v}}_{V^4(\TT^2)^2}^2 + C \pare{\norm{\widehat{w}}_{L^\infty(\TT^2)^2} + \norm{v-\widehat{v}}_{V^4(\TT^2)^2}} \norm{\Ocal}_{L^2(\TT^2)^2}^2, \notag
\end{align}
\begin{align}
	\label{eq:Ocal2D10} \abs{\psca{(w-\widehat{w})\cdot\nabla \widehat{\Vcal},\Ocal}_{L^2(\TT^2)^2}} &\leq \norm{w-\widehat{w}}_{L^\infty(\TT^2)^2} \norm{\nabla\widehat{\Vcal}}_{L^2(\TT^2)^2} \norm{\Ocal}_{L^2(\TT^2)^2}\\
	&\leq C \norm{w-\widehat{w}}_{V^3(\TT^2)^2} \norm{\widehat{v}}_{V^4(\TT^2)^2} \norm{\Ocal}_{L^2(\TT^2)^2} \notag\\
	&\leq C \norm{\widehat{v}}_{V^4(\TT^2)^2} \norm{\Ocal}_{L^2(\TT^2)^2}^2. \notag
\end{align}
For the forcing term, we have
\begin{equation}
	\label{eq:Ocal2D13} \abs{\psca{\rot(\PP f - \PP\widehat{f}),\Ocal}_{L^2(\TT^2)^2}} \leq C \pare{\norm{\PP f - \PP\widehat{f}}_{V^1(\TT^2)^2}^2 + \norm{\Ocal}_{L^2(\TT^2)^2}^2}, 
\end{equation}
and finally,
\begin{equation}
	\label{eq:Ocal2D14} \abs{\psca{\rot(w-\widehat{w}),\Ocal}_{L^2(\TT^2)^2}} \leq \norm{\rot(w-\widehat{w})}_{L^2(\TT^2)^2} \norm{\Ocal}_{L^2(\TT^2)^2} \leq C \norm{\Ocal}_{L^2(\TT^2)^2}^2.
\end{equation}

%----- %----- %----- %----- %----- %

Now, using Estimates \eqref{eq:Ocal2D01} to \eqref{eq:Ocal2D14}, we obtain
\begin{multline*}
	\frac{1}2\frac{d}{dt} \norm{\Ocal}_{L^2(\TT^2)^2}^2 + \frac{\nu}{\alpha} \norm{\Ocal}_{L^2(\TT^2)^2}^2\\ 
	\leq C \norm{\widehat{w}}_{V^4(\TT^2)^2} \norm{v-\widehat{v}}_{V^4(\TT^2)^2}^2 + C \norm{\PP f - \PP\widehat{f}}_{V^1(\TT^2)^2}^2\\
	+ C \pare{1 + \norm{\widehat{w}}_{V^4(\TT^2)^2} + \norm{\widehat{v}}_{V^4(\TT^2)^2} + \norm{v-\widehat{v}}_{V^4(\TT^2)^2}} \norm{\Ocal}_{L^2(\TT^2)^2}^2.	
\end{multline*}
For any $0 \leq t \leq T$, the Gronwall lemma implies that
\begin{multline}
	\label{eq:Ocal2DH3} \norm{\Ocal(t)}_{L^2(\TT^2)^2}^2 \leq C e^{C_2 t} \Big(\norm{w(0)-\widehat{w}(0)}_{V^3(\TT^2)^2}^2 \\ 
	+ T \norm{\PP f - \PP\widehat{f}}_{L^\infty_tV^1(\TT^2)^2}^2 + C_1 T \norm{v-\widehat{v}}_{L^\infty_tV^4(\TT^2)^2}^2\Big),
\end{multline}
where $C$ is a generic positive constant and
\begin{align*}
	&C_1 = \norm{\widehat{w}}_{L^\infty(0,T,V^4(\TT^2)^2)},\\
	&C_2 = 1 + \norm{\widehat{w}}_{L^\infty(0,T,V^4(\TT^2)^2)} + \norm{\widehat{v}}_{L^\infty(0,T,V^4(\TT^2)^2)} + \norm{v-\widehat{v}}_{L^\infty(0,T,V^4(\TT^2)^2)}.
\end{align*}
The inequality \eqref{eq:Lipschitz2Dbis} of Theorem \ref{th:FG2pert2D} is then proved.

%----- %----- %----- %----- %----- %

\bigskip

\section{Study of the extended controlled system} \label{se:Extension}

In this section, we want to show that the approximate controllability of the system \eqref{sys:FG2etaCal} is equivalent to the approximate controllability of the system \eqref{sys:FG2etazetaCal} by low-mode controls. For any finite-dimensional subspace $E$ of $V^3(\TT^2)^2$, we remark that the approximate controllability of the system \eqref{sys:FG2etaCal} by $E$-valued controls implies immediately the approximate controllability of the system \eqref{sys:FG2etazetaCal} in the same space of controls. Indeed, we only need to choose $\zeta = 0$ in the system \eqref{sys:FG2etazetaCal}. Then, in order to prove Theorem \ref{th:FG2extension}, we only need to prove the following result 
\begin{thm}
	\label{th:ExtensionP}
	Let $T>0$ and $E$ be a finite-dimensional subspace of $V^3(\TT^2)^2$. Let $\eta$, $\zeta$ in $L^\infty\pare{0,T,E}^2$ and $$\Ucal \in L^\infty\pare{0,T,V^2(\TT^2)^2} \cap C\pare{0,T,V^2(\TT^2)^2}$$ be a solution of \eqref{sys:FG2etazetaCal}. Then, for any $k \in \NN^*$, there are a control $\eta_k \in L^\infty\pare{0,T,E}^2$ and a solution $$\Ucal_k \in L^\infty\pare{0,T,V^2(\TT^2)^2} \cap C\pare{0,T,V^2(\TT^2)^2}$$ of the system \eqref{sys:FG2etaCal}, with $\eta = \eta_k$, such that $\Ucal_k(0) = \Ucal_0$, and
	$$\norm{\Ucal_k(T) - \Ucal(T)}_{V^3(\TT^2)^2} \leq \frac{1}k.$$
\end{thm}

\bigskip

\noindent\textbf{Proof. } First of all, we can rewrite the system \eqref{sys:FG2etazetaCal} as
\begin{equation}
	\label{sys:FG2etazetaP1}
	\left\{
	\begin{aligned}
		&\dd_t\Ucal + \Lcal\Ucal + \Bcal(\Ucal,\Ucal) + \Bcal(\Ucal,\zeta) + \Bcal(\zeta,\Ucal) = \PP f + \eta - \Lcal\zeta - \Bcal(\zeta,\zeta)\\
		&\divv \Ucal = 0\\
		&\Ucal(0) = \Ucal_0,
	\end{aligned}
	\right.
\end{equation}
where 
\begin{equation*}
	\Lcal\Ucal = - \nu\PP \DD (I-\alpha\DD)^{-1} \Ucal \quad \mbox{ and} \quad \Bcal(\Ucal_1,\Ucal_2) = \PP \pare{\rot \Ucal_1\times \pare{(I-\alpha\DD)^{-1} \Ucal_2}}.
\end{equation*} 
Thus, applying Theorem \ref{th:FG2pertCal} to this system, for any $\eta, \zeta \in L^{\infty}(0,T,E)^2$, we obtain the existence of a unique solution of the system \eqref{sys:FG2etazetaP1} (or \eqref{sys:FG2etazetaCal}) 
$$\overline{\Ucal} \in L^\infty\pare{0,T,V^2(\TT^2)^2} \cap C\pare{0,T,V^2(\TT^2)^2}.$$

Next, we remark that we can also rewrite the system \eqref{sys:FG2etazetaCal} as 
\begin{equation}
	\label{sys:FG2etazetaP2}
	\left\{
	\begin{aligned}
		&\dd_t\pare{\Ucal+\zeta} + \Lcal\pare{\Ucal+\zeta} + \Bcal\pare{\Ucal+\zeta} = \PP f + \tilde{\eta}\\
		&\divv \Ucal = 0\\
		&\Ucal(0) = \Ucal_0,
	\end{aligned}
	\right.
\end{equation}
where $$\tilde{\eta} = \eta + \dd_t \zeta,$$ which means that, if $\overline{\Ucal}$ is a solution of the system \eqref{sys:FG2etazetaCal} and if $\tilde{\eta}$ belongs to $L^{\infty}(0,T,E)^2$, then $\overline{\Ucal} + \zeta$ is a solution of the system \eqref{sys:FG2etaCal}, with $\eta$ replaced by $\tilde{\eta}$ and $\Ucal_0$ by $\Ucal_0 + \zeta(0)$. So, if we want to construct a solution of the controlled system \eqref{sys:FG2etaCal} satisfying the conditions of Theorem \ref{th:ExtensionP}, we only need to check whether the conditions at time $t=0$ and $t=T$ are satisfied. To this end, we will consider a sequence of controls 
$$\zeta_l \in C^1(0,T,E), \qquad \forall\; l \in \NN^*$$ such that $$\zeta_l(0) = \zeta_l(T) = 0$$ and 
$$\lim_{l\to +\infty} \norm{\zeta_l - \zeta}_{L^2(0,T,V^2(\TT^2)^2)} = 0.$$ 
Applying Theorem \ref{th:FG2pertCal} to the system \eqref{sys:FG2etazetaP1}, for any $l \in \NN^*$, there exists a unique solution 
$$\overline{\Ucal}_l \in L^\infty\pare{0,T,V^2(\TT^2)^2} \cap C\pare{0,T,V^2(\TT^2)^2}$$ 
of the system \eqref{sys:FG2etazetaP1} (or \eqref{sys:FG2etazetaCal}), with $\zeta$ replaced by $\zeta_l$. Moreover, for any $k \in \NN^*$, there exists $l_0 \in \NN^*$ such that, for any $l \geq l_0$, for any $t\in[0,T]$, we have
\begin{equation*}
	\norm{\overline{\Ucal}_l(t) - \overline{\Ucal}(t)}_{V^1} \leq C \pare{\norm{\zeta_l - \zeta}_{L^2(V^2)} + \norm{\Lcal\zeta_l - \Lcal\zeta}_{L^2(V^1)} + \norm{\Bcal(\zeta_l,\zeta_l) - \Bcal(\zeta,\zeta)}_{L^2(V^1)}} \leq \frac{1}k.
\end{equation*}

% for any $\ee > 0$, for $k_0$ large enough and for any $k \geq k_0$, 
Now, we set $$\Ucal_k = \overline{\Ucal}_{l_0} + \zeta_{l_0} \quad \mbox{ and } \eta_k = \eta + \dd_t \zeta_{l_0}.$$ Then, $\Ucal_k$ is the solution of the system \eqref{sys:FG2etaCal}, with $\eta$ replaced by $\eta_k \in L^{\infty}(0,T,E)^2$. Moreover, we have
$$\Ucal_k(0) = \overline{\Ucal}_{l_0}(0) + \zeta_{l_0}(0) = \overline{\Ucal}_{l_0}(0) = \Ucal_0,$$ and 
\begin{equation*}
	\norm{\Ucal_k(T) - \overline{\Ucal}(T)}_{V^1(\TT^2)^2} = \norm{\overline{\Ucal}_{l_0}(T) - \overline{\Ucal}(T)}_{V^1(\TT^2)^2} \leq \frac{1}k.
\end{equation*}
Theorem \ref{th:ExtensionP} is proved. \qquad $\blacksquare$

%----- %----- %----- %----- %----- %

\bigskip
\section{Convexification of the controlled system} \label{se:Convexification}

This section is devoted the to the proof of Theorem \ref{th:FG2convex}. From the definition \ref{de:Convexification} of $\Fcal(E)$, we remark that $E \subset \Fcal(E)$ and so, the approximate controllability of the system \eqref{sys:FG2eta} (or \eqref{sys:FG2etaCal}) by $E$-valued controls evidently implies the approximate controllability of the system \eqref{sys:FG2eta} (or \eqref{sys:FG2etaCal}) by $\Fcal(E)$-valued controls. In order to prove Theorem \ref{th:FG2convex}, we only need to prove that
\begin{thm}
	\label{th:ConvexP} 
	Let $T>0$, $E$ be a finite-dimensional subspace of $V^3(\TT^2)^2$. Let $\overline{\eta}$ be a control in $L^\infty\pare{0,T,\mathcal{F}(E)}^2$ and $$\overline{\Ucal} \in L^\infty\pare{0,T,V^2(\TT^2)^2} \cap C\pare{0,T,V^2(\TT^2)^2}$$ be a solution of \eqref{sys:FG2etaCal}, with $\eta$ replaced by $\overline{\eta}$. Then, for any $k \in \NN^*$, there are controls $\eta_k, \zeta_k \in L^\infty\pare{0,T,E}^2$ and a solution $$\Ucal_k \in L^\infty\pare{0,T,V^2(\TT^2)^2} \cap C\pare{0,T,V^2(\TT^2)^2}$$ of the system \eqref{sys:FG2etazetaCal}, with $\eta$ replaced by $\eta_k$ and $\zeta$ by $\zeta_k$, such that $\Ucal_k(0) = \Ucal_0$, and 
	$$\lim_{k\to +\infty} \norm{\Ucal_k(T) - \overline{\Ucal}(T)}_{V^1(\TT^2)^2} = 0.$$
\end{thm}

\bigskip

The proof of Theorem \ref{th:ConvexP} will be divided into several steps.

\bigskip

\noindent\textbf{Step 1. } \textit{Reduction of the proof of Theorem \ref{th:ConvexP} to the case of $\mathcal{F}(E)$-valued piecewise constant in time controls.}

We suppose that Theorem \ref{th:ConvexP} is true for $\mathcal{F}(E)$-valued piecewise constant controls. We want to prove that Theorem \ref{th:ConvexP} is also true in the general case. Let $\overline{\eta} \in L^\infty(0,T,\mathcal{F}(E))$ and let $$\overline{\Ucal} \in L^\infty\pare{0,T,V^2(\TT^2)^2} \cap C\pare{0,T,V^2(\TT^2)^2}$$ be a solution of \eqref{sys:FG2etaCal}, with $\eta$ replaced by $\overline{\eta}$. We consider an approximation of $\overline{\eta}$ by a sequence $\set{\eta^m}$ of $\mathcal{F}(E)$-valued piecewise constant in time controls such that
$$\lim_{m\to +\infty} \norm{\eta^m - \overline{\eta}}_{L^2(0,T,V^2(\TT^2)^2)} = 0.$$ Applying Theorem \ref{th:FG2pertCal} (while taking $\Vcal=0$ and replacing $f$ by $f+\eta^m$ and $f+\overline{\eta}$), we deduce the existence of a solution $$\Ucal^m \in L^\infty\pare{0,T,V^2(\TT^2)^2} \cap C\pare{0,T,V^2(\TT^2)^2}$$ of the system \eqref{sys:FG2etaCal}, with $\eta$ replaced by $\eta^m$, such that
\begin{equation*}
	\sup_{t\in[0,T]} \norm{\Ucal^m(t) - \overline{\Ucal}(t)}_{V^1(\TT^2)^2} \leq C \pare{\norm{\overline{\Ucal}}_{L^2(0,T,V^2(\TT^2)^2)}} \norm{\eta^m-\overline{\eta}}_{L^2(0,T,V^2(\TT^2)^2)} \leq \frac{\ee}2,
\end{equation*}
for any $m$ larger than a certain $m_0 \in \NN^*$.

If Theorem \ref{th:ConvexP} is true for $\mathcal{F}(E)$-valued piecewise constant controls, we deduce the existence of $\eta$ and $\zeta$ in $L^\infty(0,T,E)^2$, and a solution $$\Ucal \in L^\infty\pare{0,T,V^2(\TT^2)^2} \cap C\pare{0,T,V^2(\TT^2)^2}$$ of the system \eqref{sys:FG2etazetaCal} such that
$$\Ucal(0) = \Ucal^{m_0}(0) = \overline{\Ucal}(0),$$
and 
\begin{equation*}
	\sup_{t\in[0,T]} \norm{\Ucal(t) - \Ucal^{m_0}(t)}_{V^1(\TT^2)^2} \leq \frac{\ee}2.
\end{equation*}
Therefore,
\begin{equation*}
	\sup_{t\in[0,T]} \norm{\Ucal(t) - \overline{\Ucal}(t)}_{V^1(\TT^2)^2} \leq \sup_{t\in[0,T]} \norm{\Ucal(t) - \Ucal^{m_0}(t)}_{V^1(\TT^2)^2} + \sup_{t\in[0,T]} \norm{\Ucal^{m_0}(t) - \overline{\Ucal}(t)}_{V^1(\TT^2)^2} \leq \ee.
\end{equation*}

\begin{rem}
Using an argument by iteration, we can reduce the study to the case where the control $\overline{\eta}$ is constant in time. So from now on, we will consider $\overline{\eta} \in \mathcal{F}(E)$, which is constant in the time variable.
\end{rem}

%----- %----- %----- %----- %----- %

\bigskip
\noindent \textbf{Step 2. } \textit{Construction of solutions of the extended controlled system \eqref{sys:FG2etazetaCal}.}

The construction of controls $\eta_k$, $\zeta_k$ and a solution $\Ucal_k$ of the controlled system \eqref{sys:FG2etazetaCal}, with $(\eta,\zeta)$ replaced by $(\eta_k,\zeta_k)$, follows the lines of the construction in \cite{S2006} (see also \cite{AS2005} and \cite{AS2006}). We recall that 
$$\overline{\Ucal} \in L^\infty\pare{0,T,V^2(\TT^2)^2} \cap C\pare{0,T,V^2(\TT^2)^2}$$ 
is the solution of the controlled system \eqref{sys:FG2etaCal}, with $\eta$ replaced by $\overline{\eta} \in \mathcal{F}(E)$. Let $\ee > 0$ and $\delta > 0$ which will be made more precise later and choose $N > 0$ large enough such that 
$$\norm{\PP f - P_N\PP f}_{L^2(0,T,V^1(\TT^2)^2)} + \norm{\Ucal_0-P_N\Ucal_0}_{V^1(\TT^2)^2} \leq \delta,$$ where $P_N$ is the projection onto the space of the first $N$ eigenvectors of the Stokes operator $-\mathbb{P}\DD$ and where $\mathbb{P}$ is the Leray projection onto the subspace of divergence-free vector fields of $L^2(\TT^2)^2$. Let $$\Vcal_0 = P_N\Ucal_0$$ and $\Vcal_N$ be the solution of the system
\begin{equation}
	\label{sys:FG2etaN}
	\left\{ 
	\begin{aligned} 
		&\dd_t \Vcal_N + \Lcal \Vcal_N + \Bcal\pare{\Vcal_N} = \PP P_N f + \overline{\eta}\\
		&\divv \Vcal_N = 0\\
		&\Vcal_N(0) = \Vcal_0 = P_N\Ucal_0.
	\end{aligned}
	\right.
\end{equation}
Using the definition in \eqref{eq:LcalBcal}, we remark that $(I-\alpha\DD)^{-1} \Vcal_N$ is the solution of the system \eqref{sys:FG2eta}, with $u_0 = (I-\alpha\DD)^{-1} P_N \Ucal_0$ and $\eta$ replaced by $\overline{\eta}$. Then, applying [\cite{PRR}, Theorems 2.1 and 2.4], we obtain the existence of a unique solution $$\Vcal_N \in L^\infty\pare{0,T,V^3(\TT^2)^2} \cap C\pare{0,T,V^3(\TT^2)^2}$$ of the system \eqref{sys:FG2etaN}. The following lemma (see \cite{S2006}) allows us to have a ``good decomposition'' of $\mathcal{F}(E)$-valued controls in terms of $E$-valued controls.
\begin{lem}
	\label{le:newdefE1} Let $E$ be a finite-dimensional subspace of $V^3(\TT^2)^2$. Then, for any $\overline{\eta} \in \mathcal{F}(E) \setminus E$, there exist $m \in \NN^*$; $\eta, \rho^1, \ldots, \rho^m \in E$ and $\lambda_1, \ldots, \lambda_m \in \RR_+^*$, with $\sum_{j=1}^m \lambda_j = 1$, such that, for any $\Ucal\in V^2(\TT^2)^2$, we have
	\begin{equation*}
		\Bcal\pare{\Ucal} - \overline{\eta} = \sum_{j=1}^m \lambda_j \pare{\Bcal\pare{\Ucal+\rho^j} + \Lcal \rho^j} - \eta.
	\end{equation*}
\end{lem}

\noindent\textbf{Proof. } Since $\overline{\eta} \in \mathcal{F}(E) \setminus E$, Definition \ref{de:Convexification} implies that there exist $$k \in \NN^*;\; \alpha_1, \ldots, \alpha_k > 0;\; \mbox{and } \tilde{\eta}, \tilde{\rho}^1, \ldots, \tilde{\rho}^k \in E$$ such that $$\overline{\eta} = \tilde{\eta} - \sum_{j=1}^k \alpha_j \Bcal(\tilde{\rho}^j).$$ Let $m=2k$, $\alpha = \alpha_1 + \ldots + \alpha_k$ and
\begin{displaymath}
	\left\{
	\begin{aligned}
		&\lambda_j = \frac{\alpha_j}{2\alpha},\; \rho^j = \sqrt{\alpha} \tilde{\rho}^j, \quad \forall\; j\in \set{1,\ldots,k}\\
		&\lambda_j = \frac{\alpha_{j-k}}{2\alpha},\; \rho^j = -\sqrt{\alpha} \tilde{\rho}^{j-k}, \quad \forall\; j\in \set{k+1,\ldots,m}.
	\end{aligned}
	\right.
\end{displaymath}
We remark that, for any $j \in \set{1,\ldots,k}$, we have $$\lambda_j = \lambda_{j+k} \quad\mbox{ and }\quad \rho^j = - \rho^{j+k}.$$ Then, for any $\Ucal \in V^2(\TT^2)^2$, direct calculations give, $$\Bcal(\Ucal,\Ucal) - \overline{\eta} = \sum_{j=1}^m \lambda_j \pare{\Bcal(\Ucal+\rho^j) + \Lcal \rho^j} - \eta. \qquad \qquad \blacksquare$$

\noindent Lemma \ref{le:newdefE1} allows us to rewrite the system \eqref{sys:FG2etaN} as follows
\begin{equation}
	\label{sys:FG2etaN1}
	\left\{ 
	\begin{aligned}
		&\dd_t \Vcal_N + \Lcal \Vcal_N + \sum_{j=1}^m \lambda_j \pare{\Bcal(\Vcal_N+\rho^j) + \Lcal \rho^j} = \PP P_N f + \eta\\
		&\divv \Vcal_N = 0\\
		&\Vcal_N(0) = \Vcal_0 = P_N\Ucal_0.
	\end{aligned}
	\right.
\end{equation}

\bigskip

Now, we use the same construction as explained in \cite{S2006} to build the needed additional control $\zeta$ in the system \eqref{sys:FG2etazetaCal}. To this end, we introduce the following 1-periodic function $\varphi \;:\; \RR_+ \to E$
\begin{equation}
	\label{eq:phi}
	\left\{ 
	\begin{aligned}
		&\varphi(s) = \varphi(s+1) && \mbox{for any  } s\in \RR_+\\
		&\varphi(s) = \rho^1 && \mbox{si  } 0 \leq s < \lambda_1\\
		&\varphi(s) = \rho^j && \mbox{si  } \lambda_1 + \ldots + \lambda_{j-1} \leq s < \lambda_1 + \ldots + \lambda_j, \quad \mbox{for any } 2 \leq j \leq m\\
	\end{aligned}
	\right. 
\end{equation}
For any $k \geq 1$, let
\begin{equation}
	\label{eq:psik} \psi_k(t) = \varphi\pare{\frac{kt}{T}}.
\end{equation}
Then, the system \eqref{sys:FG2etaN1} can be rewritten as follows
\begin{equation}
	\label{sys:FG2etaN2}
	\left\{ 
	\begin{aligned}
		&\dd_t \Vcal_N + \Lcal (\Vcal_N + \psi_k) + \Bcal(\Vcal_N+\psi_k) = \PP P_Nf + \eta + f_k\\
		&\divv \Vcal_N = 0\\
		&\Vcal_N(0) = \Vcal_0 = P_N\Ucal_0,
	\end{aligned}
	\right.
\end{equation}
where
\begin{equation}
	\label{eq:3fk} f_k(t) = g_k(t) + h_k(t),
\end{equation}
with
\begin{equation}
	\label{eq:3ghk} 
	\left\{ 
	\begin{aligned}
		&g_k(t) = \Lcal \psi_k(t) - \sum_{j=1}^m \lambda_j \Lcal \rho^j\\
		&h_k(t) = \Bcal(\Vcal_N + \psi_k(t)) - \sum_{j=1}^m \lambda_j \Bcal(\Vcal_N + \rho^j).
	\end{aligned}
	\right. 
\end{equation}
We remark that, for any $s \in \set{1,2}$ and for any $\Ucal \in V^s(\TT^2)^2$, we have $$\norm{\Lcal\Ucal}_{V^s(\TT^2)^2} \leq \frac{\nu}{\alpha} \norm{\Ucal}_{V^s(\TT^2)^2}.$$ Then, for $s \in \set{1,2}$ and for any $t \geq 0$, simple calculations give,
\begin{align}
	\norm{g_k(t)}_{V^s(\TT^2)^2} &\leq \norm{\Lcal \psi_k(t)}_{V^s(\TT^2)^2} + \sum_{j=1}^m \lambda_j \norm{\Lcal \rho^j}_{V^s(\TT^2)^2} \leq \frac{2\nu}{\alpha} \max_{1\leq j\leq m} \norm{\rho^j}_{V^s(\TT^2)^2}\\
	\norm{h_k(t)}_{V^s(\TT^2)^2} &\leq \norm{B(\Vcal_N+\psi_k(t))}_{V^s(\TT^2)^2} + \sum_{j=1}^m \lambda_j \norm{B(\Vcal_N+\rho^j)}_{V^s(\TT^2)^2}\\ 
	&\leq 2 \max_{1\leq j\leq m} \norm{B(\Vcal_N+\rho^j)}_{V^s(\TT^2)^2}. \notag
\end{align}
Concerning the bilinear operator $\Bcal$, classical results imply, for any $\Ucal \in V^3(\TT^2)^2$ and for any $\Vcal \in V^2(\TT^2)^2$,
\begin{align}
	\label{eq:bilin01}
	\norm{\Bcal(\Ucal,\Vcal)}_{V^1(\TT^2)^2} &\leq \norm{\rot\pare{\rot \Ucal \times \pare{(I-\alpha\DD)^{-1}\Vcal}}}_{L^2(\TT^2)^2}\\
	&= \norm{\pare{(I-\alpha\DD)^{-1}\Vcal}\cdot\nabla (\rot\Ucal)}_{L^2(\TT^2)^2} \leq C(\alpha) \norm{\Vcal}_{V^2(\TT^2)^2} \norm{\Ucal}_{V^2(\TT^2)^2}  \notag
\end{align}
and
\begin{align}
	\label{eq:bilin02}
	\norm{\Bcal(\Ucal,\Vcal)}_{V^2(\TT^2)^2} &\leq C(\alpha) \Big(\norm{\rot \Ucal}_{V^2(\TT^2)^2} \norm{(I-\alpha\DD)^{-1}\Vcal}_{L^\infty(\TT^2)^2}\\
	&\qquad \qquad \qquad \qquad + \norm{\rot \Ucal}_{L^\infty(\TT^2)^2} \norm{(I-\alpha\DD)^{-1}\Vcal}_{V^2(\TT^2)^2}\Big)  \notag\\
	&\leq C(\alpha) \norm{\Vcal}_{V^2(\TT^2)^2} \norm{\Ucal}_{V^3(\TT^2)^2}.  \notag
\end{align}
Then, we obtain, for $s\in\set{1,2}$ and for any $t\geq 0$,
\begin{equation}
	\label{eq:3gk} \norm{g_k(t)}_{V^s(\TT^2)^2} \leq \frac{2\nu}{\alpha} \max_{1\leq j\leq m} \norm{\rho^j}_{V^s(\TT^2)^2}
\end{equation}
and
\begin{equation}
	\label{eq:3hk} \norm{h_k(t)}_{V^s(\TT^2)^2} \leq 2C(\alpha) \max_{1\leq j\leq m} \norm{\Vcal_N+\rho^j}_{V^2(\TT^2)^2} \norm{\Vcal_N+\rho^j}_{V^{s+1}(\TT^2)^2}.
\end{equation}

%\vspace{3cm}

%----- %----- %----- %----- %----- %

\bigskip

Next, for any $f\in L^\infty(0,T,H^1_{per}(\TT^2)^2)$, we let $\Kcal f$ be the solution of the system 
\begin{equation}
	\label{sys:FG2lin}
	\left\{ 
	\begin{aligned}
		&\dd_t \Zcal + \Lcal \Zcal = \PP f\\
		&\divv \Zcal = 0\\
		&\Zcal(0) = 0.
	\end{aligned}
	\right.
\end{equation}
By considering $z = (I-\alpha\DD)^{-1} \Zcal$ and applying [\cite{PRR}, Theorem 2.4], we have the following 
\begin{lem}
	\label{le:Kerf} Let $s \in \NN^*$. If $f\in L^\infty(0,T,V^s(\TT^2)^2)$ then $$\Kcal f \in L^\infty(0,T,V^s(\TT^2)^2) \cap C(0,T,V^s(\TT^2)^2).$$
\end{lem}
\noindent For any $k\in \NN^*$, we set
\begin{equation*}
	\Wcal_k = \Vcal_N - \Kcal f_k.
\end{equation*}
Then, the system \eqref{sys:FG2etaN2} becomes
\begin{equation}
	\label{sys:WKF}
	\left\{ 
	\begin{aligned}
		&\dd_t \Wcal_k + \Lcal\Wcal_k + \Bcal(\Wcal_k) + \Bcal(\Wcal_k,\psi_k + \Kcal f_k) + \Bcal(\psi_k + \Kcal f_k,\Wcal_k)\\ 
		&\qquad \qquad \qquad \qquad \qquad \qquad \qquad \qquad \qquad \qquad \qquad = \PP P_Nf + \eta - \PP B(\psi_k + \Kcal f_k) - \Lcal\psi_k\\
		&\divv \Wcal_k = 0\\
		&\Wcal_k(0) = \Vcal_N(0) = P_N\Ucal_0.
	\end{aligned}
	\right. 
\end{equation}
In other words, $\Wcal_k$ is the solution of the system \eqref{sys:FG2perCal}, with data $(\Vcal,f,\Wcal_0)$ replaced by
\begin{equation*}
	(\psi_k + \Kcal f_k, P_Nf + \eta - \Bcal(\psi_k + \Kcal f_k) - \Lcal \psi_k, P_N\Ucal_0).
\end{equation*}
Let
\begin{equation*}
	\overline{\Vcal} = \psi_k, \quad \overline{f} = f + \eta - \Bcal(\psi_k) - \Lcal \psi_k, \quad \overline{\Wcal}_0 = \Ucal_0,
\end{equation*}
and let $\Ucal_k$ be the solution of the system \eqref{sys:FG2perCal} with data $(\overline{\Vcal},\overline{f},\overline{\Wcal}_0)$. It is then easy to show that $\Ucal_k$ is the solution of the controlled system \eqref{sys:FG2etazetaCal}, with controls $\eta$ and $\zeta = \psi_k$
\begin{equation*}
	\left\{
	\begin{aligned}
		&\dd_t \Ucal_k + \Lcal \pare{\Ucal_k+\psi_k} + \Bcal\pare{\Ucal_k+\psi_k} = \PP f + \eta\\
		&\divv \Ucal_k = 0\\
		&\Ucal_k(0) = \Ucal_0.
	\end{aligned}
	\right.
\end{equation*}
So, all we need to do now is to prove that, for any $\ee > 0$, there exists $k_0 \in \NN^*$ such that, for any $k \geq k_0$, we have
$$\norm{\Ucal_k(T) - \overline{\Ucal}(T)}_{V^1(\TT^2)^2} \leq \ee,$$ where $\overline{\Ucal}$ is the solution of the controlled system \eqref{sys:FG2etaCal}, with $\eta$ replaced by $\overline{\eta} \in \mathcal{F}(E)$.

\bigskip

\noindent Recall that $\Wcal_k = \Vcal_N - \Kcal f_k$. For any $t \geq 0$, we have
\begin{align}
	\label{eq:UKUBar} \norm{\Ucal_k(t) - \overline{\Ucal}(t)}_{V^1(\TT^2)^2} &\leq \norm{\Ucal_k(t) - \Vcal_N(t)}_{V^1(\TT^2)^2} + \norm{\Vcal_N(t) - \overline{\Ucal}(t)}_{V^1(\TT^2)^2}\\
	&\leq \norm{\Ucal_k(t) - \Wcal_k(t)}_{V^1(\TT^2)^2} + \norm{\Kcal f_k(t)}_{V^1(\TT^2)^2} + \norm{\Vcal_N(t) - \overline{\Ucal}(t)}_{V^1(\TT^2)^2}. \notag
\end{align}
Applying Theorem \ref{th:FG2pertCal}, we obtain
\begin{equation*}
	\norm{\Vcal_N(t) - \overline{\Ucal}(t)}_{V^1(\TT^2)^2} \leq C \pare{\norm{f - P_Nf}_{L^2(0,T,V^1(\TT^2)^2)} + \norm{\Ucal_0-P_N\Ucal_0}_{V^1(\TT^2)^2}} \leq C\delta.
\end{equation*}
%Applying Theorem \ref{th:FG2pertCal} once again with
% \begin{equation*}
% 	\widehat{\Vcal} = \psi_k + \Kcal f_k, \quad \widehat{f} = P_N f + \eta - \Bcal(\psi_k + \Kcal f_k) - \Lcal\psi_k, \quad \widehat{\Wcal}_0 = P_N\Ucal_0.
% \end{equation*}
% and
% \begin{equation*}
% 	\overline{\Vcal} = \psi_k, \quad \overline{f} = f + \eta - \Bcal(\psi_k) - \Lcal\psi_k, \quad \overline{\Wcal}_0 = \Ucal_0,
% \end{equation*}
% we get
and
\begin{align*}
	&\norm{\Ucal_k(t) - \Wcal_k(t)}_{V^1(\TT^2)^2} \leq C\Big(\norm{\Kcal f_k}_{L^2(0,T,V^2(\TT^2)^2)} + \norm{\Bcal(\psi_k + \Kcal f_k) - \Bcal(\psi_k)}_{L^2(0,T,V^1(\TT^2)^2)} \\
	&\qquad \qquad \qquad \qquad \qquad \qquad \qquad + \norm{f - P_Nf}_{L^2(0,T,V^1(\TT^2)^2)} + \norm{\Ucal_0 - P_N\Ucal_0}_{V^1(\TT^2)^2}\Big). \notag
\end{align*}
Next, Estimate \eqref{eq:bilin01} implies that
\begin{align*}
	&\norm{\Bcal(\psi_k + \Kcal f_k) - \Bcal(\psi_k)}_{L^2(0,T,V^1(\TT^2)^2)}\\ 
	&\qquad \qquad \leq \norm{\Bcal(\psi_k, \Kcal f_k)}_{L^2(0,T,V^1(\TT^2)^2)} + \norm{\Bcal(\Kcal f_k,\psi_k)}_{L^2(0,T,V^1(\TT^2)^2)} + \norm{\Bcal(\Kcal f_k)}_{L^2(0,T,V^1(\TT^2)^2)}\\
	&\qquad \qquad \leq CT^{\frac{1}2} \norm{\Kcal f_k}_{L^\infty(0,T,V^2(\TT^2)^2)} \pare{\norm{\Kcal f_k}_{L^\infty(0,T,V^2(\TT^2)^2)} + \max_{1\leq j\leq m} \norm{\rho^j}_{V^2(\TT^2)^2}}
\end{align*}
Then, we deduce from \eqref{eq:UKUBar} that
\begin{align*}
	&\norm{\Ucal_k(t) - \overline{\Ucal}(t)}_{V^1(\TT^2)^2}\\ 
	&\qquad \qquad \leq C\delta + CT^{\frac{1}2} \norm{\Kcal f_k}_{L^\infty(0,T,V^2(\TT^2)^2)} \pare{\norm{\Kcal f_k}_{L^\infty(0,T,V^2(\TT^2)^2)} + \max_{1\leq j\leq m} \norm{\rho^j}_{V^2(\TT^2)^2} + 1}.
\end{align*}

\bigskip

Now, we choose $\delta > 0$ such that $C\delta \leq \frac{\ee}2$. In order to prove Theorem \ref{th:ConvexP}, we only need to prove the following lemma.
\begin{lem}
	\label{le:Kernel} We have $$\lim_{k\to +\infty} \norm{\Kcal f_k}_{L^\infty(0,T,V^2(\TT^2)^2)} = 0.$$
\end{lem}

\noindent In order to prove Lemma \ref{le:Kernel}, we need to prove the following result for $f_k$.
\begin{lem}
	\label{le:fkrelax} Let $T > 0$ and for any $k\in\NN^*$, let $f_k$ be defined as in \eqref{eq:3fk}. Then, we have
	$$\lim_{k\to +\infty} \sup_{t\in[0,T]} \norm{\int_0^t f_k(s) ds}_{V^2(\TT^2)^2} = 0.$$
\end{lem}

%----- %----- %----- %----- %----- %

\bigskip
\noindent \textbf{Step 3. } \textit{Proof of Lemma \ref{le:fkrelax}.}

For any $k\in \NN^*$ and for any $t > 0$, let $$F_k(t) = \int_0^t f_k(s) ds.$$ So, our goal is to prove that
\begin{equation}
	\label{eq:fkprim} \lim_{k\to +\infty} \norm{F_k}_{C(0,T,V^2(\TT^2)^2)} = 0. 
\end{equation}
We remark that Estimates \eqref{eq:bilin02}, \eqref{eq:3gk}, \eqref{eq:3hk} and the definition of $f_k$ imply that if \eqref{eq:fkprim} is true for all piecewise constant (with respect to the time variable) functions $\Vcal_N$, then \eqref{eq:fkprim} is true for all functions $\Vcal_N$ (by using an approximation of $\Vcal_N$ by piecewise constant functions). For this reason, we suppose that there exist $L \in \NN^*$, $t_0, \ldots, t_L \in \RR_+$ such that $$0 = t_0 < \ldots < t_L = T,$$ and that $$\Vcal_N(t) = v_q, \qquad \forall\; t\in ]t_{q-1},t_q[, \quad \forall\; q\in \set{1,\ldots,L}.$$

%----- %----- %----- %----- %----- %

\vspace{0.3cm}

Now, by using direct calculations, we can prove that $$\lim_{k\to +\infty} F_k(t) = 0.$$ For the details, we send the reader to the book of Jurdjevic \cite{Ju}. 
We remark that, for any $t \in [0,T]$, the set $\set{F_k(t)}_k$ is relatively compact in $V^2(\TT^2)^2$. Indeed, the set $\set{f_k(t)}$ only takes a finite number of values, independently of $k$. Let $M$ be the set of value of $\set{f_k(t)}$ and we suppose that $$M = \set{M_1,\ldots,M_K}, \qquad K\in \NN^*.$$ Then, there exist positive constant $a_1$, \ldots, $a_K$ such that $$a_1 + \ldots + a_K = t$$ and $$F_k(t) = \sum_{i=1}^K a_i M_i.$$ Thus, $\set{F_k(t)}_k$ is relatively compact in $V^2(\TT^2)^2$. Moreover, from \eqref{eq:3gk} and \eqref{eq:3hk}, there exists a positive constant $C_0$ such that
$$\sup_{t\in[0,T]} \norm{f_k(t)}_{V^2(\TT^2)^2} \leq C_0,$$ which means that $\set{F_k(\cdot)}_k$ is equicontinuous. Then, the Ascoli's theorem implies that $\set{F_k(\cdot)}_k$ is relatively compact in $C(0,T,V^2(\TT^2)^2)$.

%----- %----- %----- %----- %----- %

\vspace{0.3cm}

Next, we have
\begin{itemize}
	\item $\set{F_k(\cdot)}_k$ is relatively compact in $C(0,T,V^2(\TT^2)^2)$.
	\item $F_k(t) \to 0$ in $V^2(\TT^2)^2$, for any $t\in[0,T]$.
\end{itemize}
Then, it is clear that
$$\lim_{k\to +\infty} \sup_{t\in[0,T]} \norm{F_k(t)}_{V^2(\TT^2)^2} = \lim_{k\to +\infty} \sup_{t\in[0,T]} \norm{\int_0^t f_k(s) ds}_{V^2(\TT^2)^2} = 0. \qquad \blacksquare$$

%----- %----- %----- %----- %----- %

\bigskip
\noindent \textbf{Step 4. } \textit{Proof of Lemma \ref{le:Kernel}.}

We recall that $\Zcal=\Kcal f_k$ is the solution of the system
\begin{equation}
	\label{sys:FG2linKCal}
	\left\{ 
	\begin{aligned}
		&\dd_t \Zcal + \Lcal \Zcal = \PP f_k\\
		&\divv \Zcal = 0\\
		&\Zcal(0) = 0.
	\end{aligned}
	\right.
\end{equation}
or equivalently, $z = (I - \alpha \DD)^{-1} \Zcal$ is solution of the system
\begin{equation}
	\label{sys:FG2linK}
	\left\{ 
	\begin{aligned}
		&\dd_t (z-\alpha\DD z) - \nu\DD z = \PP f_k\\
		&\divv z = 0\\
		&z(0) = 0.
	\end{aligned}
	\right.
\end{equation}
We also recall that \emph{a priori} estimates in this paragraph can be justified by applying an approximation by a Galerkin scheme. Following the method presented in \cite{PRR}, we apply the $rot$ operator to the first equation of \eqref{sys:FG2linK} and then we take the $L^2$ scalar product of the obtained equation with $-\rot\!\pare{\DD z -\alpha\DD^2 z}$. We get
\begin{multline}
	\label{eq:SKH4a} \frac{1}2 \frac{d}{dt} \norm{\nabla\pare{\rot z - \alpha\;\rot\DD z}}_{L^2(\TT^2)^2}^2 + \nu\pare{\norm{\DD\rot z}_{L^2(\TT^2)^2}^2 + \alpha\norm{\DD^2 z}_{L^2(\TT^2)^2}^2}\\ 
	= -\psca{\rot f_k \;,\; \rot\!\pare{\DD z -\alpha\DD^2 z}}_{L^2(\TT^2)^2}.
\end{multline}
Integrating over $[0,t]$ and then, performing multiple integrations by parts (with respect to the space variable $x$ and then with respect to the time variable $t$), we have
\begin{align}
	\label{eq:SKH4b}
	&\frac{1}2\norm{\nabla\pare{\rot z(t) - \alpha\;\rot\DD z(t)}}_{L^2(\TT^2)^2}^2 + \nu\pare{\int_0^t \norm{\DD\rot z(s)}_{L^2(\TT^2)^2}^2 ds + \alpha\int_0^t \norm{\DD^2 z(s)}_{L^2(\TT^2)^2}^2 ds}\\
	&\qquad \qquad = - \int_0^t \int_{\TT^2} \rot f_k(s,x) \cdot \rot (\DD z - \alpha\DD^2 z)(s,x) dx ds \notag\\
	&\qquad \qquad = \int_0^t \int_{\TT^2} \nabla\rot\! f_k(s,x) : \nabla (\rot z -\alpha\; \rot \DD z)(s,x) dx dt \notag\\
	&\qquad \qquad = \int_{\TT^2} \pare{\int_0^t \nabla\rot\! f_k(s,x) ds} : \nabla \pare{\rot z - \alpha\;\rot\DD z}(t,x) dx \notag\\
	&\qquad \qquad \qquad \qquad - \int_{\TT^2} \int_0^t \pare{\int_0^s \nabla\rot\! f_k(\tau,x) d\tau} : \pare{\frac{\dd}{\dd s} \nabla\pare{\rot z - \alpha\;\rot\DD z}(s,x)} ds dx \notag\\
	&\qquad \qquad = J_1(t,k) + J_2(t,k). \notag
\end{align}

\noindent For the first term on the right-hand side, Cauchy-Schwarz inequality implies that
\begin{align}
	\label{eq:J1}
	\abs{J_1(t,k)} &\leq \norm{\int_0^t\nabla\rot\! f_k(s,x) ds}_{L^2(\TT^2)^2} \norm{\nabla\pare{\rot z(t) - \alpha\;\rot\DD z(t)}}_{L^2(\TT^2)^2}\\
	&\leq C \norm{\int_0^t f_k(s,x) ds}_{V^2(\TT^2)^2} \norm{\nabla\pare{\rot z(t) - \alpha\;\rot\DD z(t)}}_{L^2(\TT^2)^2}. \notag
\end{align}

\noindent For the second term, also using Cauchy-Schwarz inequality, we obtain
\begin{align*}
	\abs{J_2(t,k)} &\leq \int_0^t \pare{\norm{\int_0^s \nabla \rot\! f_k(\tau) d\tau}_{L^2(\TT^2)^2} \norm{\frac{\dd}{\dd s} \nabla \pare{\rot z - \alpha\;\rot\DD z}(s)}_{L^2(\TT^2)^2}} ds\\
	&\leq \pare{\int_0^t \norm{\int_0^s f_k(\tau) d\tau}_{V^2(\TT^2)^2}^2 ds}^{\frac{1}2} \pare{\int_0^t \norm{\frac{\dd}{\dd s} \nabla \pare{\rot z - \alpha\;\rot\DD z}(s)}_{L^2(\TT^2)^2}^2 ds}^{\frac{1}2}\\
	&\leq T^{\frac{1}2} \sup_{t\in[0,T]} \norm{\int_0^t f_k(s) ds}_{V^2(\TT^2)^2} \pare{\int_0^t \norm{\frac{\dd}{\dd s} \nabla \pare{\rot z - \alpha\;\rot\DD z}(s)}_{L^2(\TT^2)^2}^2 ds}^{\frac{1}2}.
\end{align*}

\noindent Now, we come back the first equation of \eqref{sys:FG2linK}. Applying the $rot$ operator to this equation and then, taking the $L^2$ scalar product of the obtained equation with $-\frac{\dd}{\dd t}\rot\!\pare{\DD z -\alpha\DD^2 z}$, we get
\begin{multline*}
	\norm{\dd_t\nabla\rot\!\pare{z-\alpha\DD z}}_{L^2(\TT^2)^2}^2 + \frac{\nu}2 \frac{d}{dt} \pare{\norm{\rot\DD z}_{L^2(\TT^2)^2}^2 + \alpha \norm{\DD^2 z}_{L^2(\TT^2)^2}^2}\\
	= \psca{\nabla\rot f_k \;,\; \dd_t\nabla\rot\!\pare{z-\alpha\DD z}}_{L^2(\TT^2)^2}. \qquad
\end{multline*}
Thus,
\begin{align*}
	&\int_0^t \norm{\dd_s\nabla\rot\!\pare{z(s)-\alpha\DD z(s)}}_{L^2(\TT^2)^2}^2 ds + \frac{\nu}2 \pare{\norm{\rot\DD z(t)}_{L^2(\TT^2)^2}^2 + \alpha \norm{\DD^2 z(t)}_{L^2(\TT^2)^2}^2}\\
	& \qquad\qquad\qquad \leq \int_0^t\norm{\nabla\rot f_k(s)}_{L^2(\TT^2)^2} \norm{\dd_s\nabla\rot\!\pare{z(s)-\alpha\DD z(s)}}_{L^2(\TT^2)^2} ds\\
	& \qquad\qquad\qquad \leq \norm{f_k}_{L^2(0,T,V^2(\TT^2)^2)} \pare{\int_0^t \norm{\dd_s\nabla\rot\!\pare{z(s)-\alpha\DD z(s)}}_{L^2(\TT^2)^2}^2 ds}^{\frac{1}2}
\end{align*}
Come back to $J_2(t,k)$, we have
\begin{equation}
	\label{eq:J2} \abs{J_2(t,k)} \leq C T \sup_{t\in[0,T]} \norm{\int_0^t f_k(s) ds}_{V^2(\TT^2)^2} \norm{f_k}_{L^\infty(0,T,V^2(\TT^2)^2)}.
\end{equation}

\noindent Combining \eqref{eq:SKH4b} with \eqref{eq:J1}, \eqref{eq:J2} and the fact that $$\lim_{k\to +\infty} \sup_{t\in[0,T]} \norm{\int_0^t f_k(s) ds}_{V^2(\TT^2)^2} = 0,$$
we conclude that
\begin{equation*}
	0 \leq \lim_{k\to +\infty} \norm{\Kcal f_k}_{L^\infty(0,T,V^2(\TT^2)^2)} \leq C \lim_{k\to +\infty} \norm{\nabla\pare{\rot z(t) - \alpha\;\rot\DD z(t)}}_{L^\infty(0,T,L^2(\TT^2)^2)} = 0. \qquad \blacksquare
\end{equation*}

%----- %----- %----- %----- %----- %

\bigskip

\section{Saturation property for the controlled system of fluids of second grade} \label{se:Saturation}

Let $q = \pare{q_1,q_2} \in \; ]0,+\infty[^2$ be fixed and 
$$\TT^2_q = \RR^2 / \ZZ^2_q \qquad \mbox{ with } \quad \ZZ^2_q = \set{x = \pare{x_1,x_2}\in \RR^2 \;\Big\vert\; \frac{x_i}{q_i} \in \ZZ, \; i = 1,2}.$$ For any $x=(x_1,x_2), y=(y_1,y_2) \in \RR^2$, let
$$\psca{x,y} = x_1y_1 + x_2y_2, \quad \abs{x} = \abs{x_1} + \abs{x_2}, \quad \norm{x} = \sqrt{\psca{x,x}},$$
and
$$\psca{x,y}_q = \frac{x_1y_1}{q_1} + \frac{x_2y_2}{q_2}, \quad \norm{x}_q = \sqrt{\frac{x_1^2}{q_1^2} + \frac{x_2^2}{q_2^2}}.$$
For any $a = (a_1,a_2) \in \RR^2 \setminus \set{0}$, let $a^\perp = (-a_2,a_1)$. We will denote $a^{q,\perp}$ the unit vector which satisfied $\psca{a,a^{q,\perp}}_q = 0$ and $\norm{a^{q,\perp}} = 1$ and we denote $P_a$ the orthogonal of $\RR^2$ onto the subspace $Span\set{a^{q,\perp}}$ generated by $a^{q,\perp}$. Direct calculations also give 
\begin{lem}
	\label{le:SaturaProj} Let $a, l \in \RR^2 \setminus \set{0}$. Then,
	\begin{align}
		\label{eq:SaturaCos} &\mathbb{P} \pare{a\cos\psca{l,x}_q} = \pare{P_l a} \cos\psca{l,x}_q\\
		\label{eq:SaturaSin} &\mathbb{P} \pare{a\sin\psca{l,x}_q} = \pare{P_l a} \sin\psca{l,x}_q
	\end{align}
\end{lem}

We recall that for any $m \in \ZZ^2 \setminus \set{0}$, we set $$c_m(x) = m^{q,\perp}\cos\psca{m,x}_q \qquad \mbox{ and } \qquad s_m(x) = m^{q,\perp}\sin\psca{m,x}_q.$$ These vector fields $c_m$, $s_m$, with $m \in \ZZ^2 \setminus \set{0}$, are eigenvectors of the Stokes operator $-\mathbb{P}\DD$ and the family $\set{c_m,s_m \;\vert\; m \in \ZZ^2 \setminus \set{0}}$ forms an orthonormal basis of $V^k(\TT^2_q)^2$, $k\in \NN$. In Section \ref{se:Intro}, for any $N\in\NN^*$, we already set 
\begin{equation}
	\tag{\ref{eq:HN}} \Hcal^N_q = Span\set{c_m,s_m \;\vert\; m\in \ZZ\setminus\set{0}, \abs{m} \leq N}.
\end{equation}

\noindent For any $m \in \ZZ^2 \setminus \set{0}$, let $$\mathcal{C}_m = Span\set{c_m,c_{-m}} \qquad\mbox{ and }\qquad \mathcal{S}_m = Span\set{s_m,s_{-m}}.$$ 
\begin{lem}
	\label{le:SaturaDecomp} Let $m, n \in \ZZ^2 \setminus \set{0}$. For any $f_m \in \mathcal{C}_m$ and $g_n \in \mathcal{S}_n$, there exist $\tilde{f}_m, \tilde{g}_n \in \RR^2$ such that
	\begin{gather*}
		\psca{\tilde{f}_m,m}_q = \psca{\tilde{g}_n,n}_q = 0,\\
		f_m(x) = \tilde{f}_m \cos\psca{m,x}_q,\\
		g_n(x) = \tilde{g}_n \sin\psca{n,x}_q.
	\end{gather*}
\end{lem}

\noindent For any $m, n \in \ZZ^2 \setminus \set{0}$ and for any $f_m \in \mathcal{C}_m$ and $g_n \in \mathcal{S}_n$, Lemma \ref{le:SaturaDecomp} allows to calculate
\begin{align*}
	&\mathbb{P}\pare{\rot f_m \times \pare{(I-\alpha\DD)^{-1} g_n}}\\ 
	&\qquad \qquad = \mathbb{P}\set{\pare{\nabla\times \pare{\tilde{f}_m\cos\psca{m,x}_q}} \times \pare{(Id-\alpha\DD)^{-1} \pare{\tilde{g}_n\sin\psca{n,x}_q}}}\\
	&\qquad \qquad = \mathbb{P}\set{\pare{1+\alpha\norm{n}_q^2}^{-1}\pint{\nabla\times\pare{\tilde{f}_m\cos\psca{m,x}_q}} \times \pare{\tilde{g}_n\sin\psca{n,x}_q}}\\
	&\qquad \qquad = \pare{1+\alpha\norm{n}_q^2}^{-1} \mathbb{P}\set{\pare{\psca{\tilde{f}_m^\perp,m}_q \sin\psca{m,x}_q} \times \pare{\tilde{g}_n\sin\psca{n,x}_q}}\\
	&\qquad \qquad = \pare{1+\alpha\norm{n}_q^2}^{-1} \psca{\tilde{f}_m^\perp,m}_q \mathbb{P} \pare{\tilde{g}_n^\perp \sin\psca{m,x}_q \sin\psca{n,x}_q}\\
	&\qquad \qquad = \frac{\pare{1+\alpha\norm{n}_q^2}^{-1}}{2} \psca{\tilde{f}_m^\perp,m}_q \mathbb{P} \pint{\tilde{g}_n^\perp \pare{\cos\psca{m-n,x}_q - \cos\psca{m+n,x}_q}}.
\end{align*}
Using Lemma \ref{le:SaturaProj}, we obtain
\begin{multline}
	\label{eq:fg} \Bcal(f_m,g_n) = \mathbb{P}\pare{\rot f_m \times \pare{(I-\alpha\DD)^{-1} g_n}}\\ = \frac{\pare{1+\alpha\norm{n}_q^2}^{-1}}{2} \psca{\tilde{f}_m^\perp,m}_q \pare{\cos\psca{m-n,x}_q P_{m-n} - \cos\psca{m+n,x}_q P_{m+n}} \tilde{g}_n^\perp.
\end{multline}
Similar calculations give
\begin{equation}
	\label{eq:gf} \Bcal(g_m,f_n) = -\frac{\pare{1+\alpha\norm{n}_q^2}^{-1}}{2} \psca{\tilde{g}_m^\perp,m}_q \pare{\cos\psca{m+n,x}_q P_{m+n} + \cos\psca{m-n,x}_q P_{m-n}} \tilde{f}_n^\perp,
\end{equation}
\begin{equation}
	\label{eq:ff} \Bcal(f_m,f_n) = \frac{\pare{1+\alpha\norm{n}_q^2}^{-1}}{2} \psca{\tilde{f}_m^\perp,m}_q \pare{\sin\psca{m+n,x}_q P_{m+n} + \sin\psca{m-n,x}_q P_{m-n}} \tilde{f}_n^\perp,
\end{equation}
and
\begin{equation}
	\label{eq:gg} \Bcal(g_m,g_n) = \frac{\pare{1+\alpha\norm{n}_q^2}^{-1}}{2} \psca{\tilde{g}_m^\perp,m}_q \pare{\sin\psca{m-n,x}_q P_{m-n} - \sin\psca{m+n,x}_q P_{m+n}} \tilde{g}_n^\perp.
\end{equation}

The next lemma is the most important result of this section, which allows us to prove the saturation property of the space of controls (Theorem \ref{th:FG2Saturation}).
\begin{lem}
	\label{le:Saturation} Let $q = (q_1,q_2)$, $q_1, q_2 > 0$. For any $m, n \in \ZZ^2 \setminus \set{0}$ satisfying
	\begin{itemize}
		\item $\ds \norm{m}_q \neq \norm{n}_q$,
		\item $m$, $n$ are not parallel,
	\end{itemize}
	and for any $f \in \mathcal{C}_{m+n}$, $g \in \mathcal{S}_{m+n}$, there exist $$a, b \in Span\set{\mathcal{C}_m, \mathcal{C}_n, \mathcal{S}_m, \mathcal{S}_n}$$ such that
	\begin{equation}
		\label{eq:Saturation00} \Bcal(a) + f, \; \Bcal(b) + g \in Span\set{\mathcal{C}_{m-n}, \mathcal{S}_{m-n}},
	\end{equation}
	where $\Bcal$ is defined in \eqref{eq:LcalBcal}.
\end{lem}

\noindent\textbf{Proof. } Taking $m=n$ in Estimates \eqref{eq:ff} and \eqref{eq:gg}, we have
\begin{align}
	\label{eq:ffbis} \Bcal\pare{f_m,f_m} &= \frac{1}{2} \pare{1+\alpha\norm{m}_q^2}^{-1} \psca{\tilde{f}_m^\perp,m}_q \sin\psca{2m,x}_q P_{2m} \tilde{f}_m^\perp,\\
	\label{eq:ggbis} \Bcal\pare{g_m,g_m} &= -\frac{1}{2} \pare{1+\alpha\norm{m}_q^2}^{-1} \psca{\tilde{g}_m^\perp,m}_q \sin\psca{2m,x}_q P_{2m} \tilde{g}_m^\perp.
\end{align}
Since $\psca{\tilde{f}_m,m}_q = \psca{\tilde{g}_m,m}_q = 0$, we deduce that $\tilde{f}_m, \tilde{g}_m \in Span\set{m_q^\perp}$. By definition, $P_{2m}$ is the projection onto the subspace $Span\set{\pare{2m}_q^\perp} = Span\set{m_q^\perp}$. So, we have $$P_{2m}\pare{\tilde{f}_m^\perp} = P_{2m}\pare{\tilde{g}_m^\perp} = 0,$$ which means that, for any $m\in\ZZ^2\setminus\set{0}$,
\begin{equation}
	\label{eq:BffBgg} \Bcal(f_m,f_m) = \Bcal(g_m,g_m) = 0.
\end{equation}

Now, following the idea of \cite{S2006}, for any $f \in \mathcal{C}_{m+n}$, we look for $a\in Span\set{\Ccal_m,\Scal_n}$ under the form 
\begin{equation}
	\label{eq:formA} a = f_m + g_n, \quad f_m \in \Ccal_m, \quad g_n \in \Scal_n,
\end{equation}
such that
$$\Bcal(a) + f \in Span\set{\mathcal{C}_{m-n}, \mathcal{S}_{m-n}}.$$
Since $a = f_m + g_n$, taking into account \eqref{eq:fg}, \eqref{eq:gf} and \eqref{eq:BffBgg}, we have
\begin{align*}
	\Bcal(a) &= \frac{1}2 \cos\psca{m-n,x}_q P_{m-n} \set{\pare{1+\alpha\norm{n}_q^2}^{-1} \psca{\tilde{f}_m^\perp,m}_q \tilde{g}_n^\perp - \pare{1+\alpha\norm{m}_q^2}^{-1} \psca{\tilde{g}_n^\perp,n}_q \tilde{f}_m^\perp}\\
	&- \frac{1}2 \cos\psca{m+n,x}_q P_{m+n} \set{\pare{1+\alpha\norm{n}_q^2}^{-1} \psca{\tilde{f}_m^\perp,m}_q \tilde{g}_n^\perp + \pare{1+\alpha\norm{m}_q^2}^{-1} \psca{\tilde{g}_n^\perp,n}_q \tilde{f}_m^\perp}.
\end{align*}
We remark that $$\cos\psca{m-n,x}_q P_{m-n} \set{\pare{1+\alpha\norm{n}_q^2}^{-1} \psca{\tilde{f}_m^\perp,m}_q \tilde{g}_n^\perp - \pare{1+\alpha\norm{m}_q^2}^{-1} \psca{\tilde{g}_n^\perp,n}_q \tilde{f}_m^\perp}$$ belongs to $Span\set{\mathcal{C}_{m-n}, \mathcal{S}_{m-n}}$. So, we only need to prove that, for any $f \in \Ccal_{m+n}$, there are $f_m\in \Ccal_m$ and $g_n \in \Scal_n$ such that
\begin{equation}
	\label{eq:ProjF} f = \frac{1}2 \cos\psca{m+n,x}_q P_{m+n} \set{\pare{1+\alpha\norm{n}_q^2}^{-1} \psca{\tilde{f}_m^\perp,m}_q \tilde{g}_n^\perp + \pare{1+\alpha\norm{m}_q^2}^{-1} \psca{\tilde{g}_n^\perp,n}_q \tilde{f}_m^\perp}.
\end{equation}
Let 
$$F = \pare{1+\alpha\norm{n}_q^2}^{-1} \psca{\tilde{f}_m^\perp,m}_q \tilde{g}_n^\perp + \pare{1+\alpha\norm{m}_q^2}^{-1} \psca{\tilde{g}_n^\perp,n}_q \tilde{f}_m^\perp.$$
Then, we will prove \eqref{eq:ProjF} if we can find $f_m\in \Ccal_m$ and $g_n \in \Scal_n$ such that $P_{m+n} F \neq 0$, or equivalently, $\psca{F,G} \neq 0$, for some vector $G \neq 0$ and $G \in Span\set{(m+n)_q^\perp}$. For the sake of simplicity, we choose
$$G = \pare{(m_2+n_2)q_1,-(m_1+n_1)q_2}.$$ 
For any $f_m\in \Ccal_m$ and $g_n \in \Scal_n$, recall that $\psca{\tilde{f}_m,m}_q = \psca{\tilde{g}_n,n}_q = 0$. Thus, there exist $C_f, C_g \in \RR$ such that
\begin{equation*}
	\tilde{f}_m = C_f \pare{m_2q_1,-m_1q_2} \qquad\mbox{and}\qquad \tilde{g}_n = C_g \pare{n_2q_1,-n_1q_2}
\end{equation*}
We get
$$\tilde{f}_m^\perp = C_f \pare{m_1q_2,m_2q_1} \qquad\mbox{and}\qquad \tilde{g}_n^\perp = C_g \pare{n_1q_2,n_2q_1},$$
and
\begin{align*}
	\psca{\tilde{f}_m^\perp,m}_q &= \frac{q_2}{q_1} m_1^2 + \frac{q_1}{q_2} m_2^2 = q_1q_2 \norm{m}_q^2\\
	\psca{\tilde{g}_n^\perp,n}_q &= \frac{q_2}{q_1} n_1^2 + \frac{q_1}{q_2} n_2^2 = q_1q_2 \norm{n}_q^2.
\end{align*}
Let 
$$C = \pare{1+\alpha\norm{m}_q^2} \pare{1+\alpha\norm{n}_q^2}.$$
Then,
\begin{equation*}
	F = C C_f C_g q_1 q_2 \pint{\pare{1+\alpha\norm{m}_q^2}\norm{m}_q^2 \pare{n_1q_2,n_2q_1} + \pare{1+\alpha\norm{n}_q^2}\norm{n}_q^2 \pare{m_1q_2,m_2q_1}},
\end{equation*}
Let $$M_q = \pare{1+\alpha\norm{m}_q^2}\norm{m}_q^2 \qquad\mbox{and}\qquad N_q = \pare{1+\alpha\norm{n}_q^2}\norm{n}_q^2.$$ Since $\norm{m}_q \neq \norm{n}_q$ and $m$, $n$ are not parallel in $\RR^2$, we finally obtain 
\begin{align*}
	\psca{F,G} &= C C_f C_g q_1 q_2 \pint{\pare{M_qn_1q_2+N_qm_1q_2} \pare{m_2q_1+n_2q_1} - \pare{M_qn_2q_1+N_qm_2q_1} \pare{m_1q_2+n_1q_2}}\\
	&= C C_f C_g q_1^2 q_2^2 \pare{M_q-N_q} \pare{n_1m_2-n_2m_1} \neq 0,
\end{align*}
if $f_m$ and $g_n$ are not zero.

To prove the second part of the lemma, for any $g \in \mathcal{S}_{m+n}$, we can look for $b$ under the form
$$b = f_m + f_n \qquad\mbox{or}\qquad b = g_m + g_n,$$ where $f_m \in \Ccal_m$, $f_n \in \Ccal_n$, $g_m \in \Scal_m$ and $g_n \in \Scal_n$. Indeed, taking into account \eqref{eq:ff}, \eqref{eq:gg} and \eqref{eq:BffBgg}, we have
\begin{align*}
	&\Bcal(f_m+f_n)\\ 
	&\quad = \frac{1}2 \sin\psca{m-n,x}_q P_{m-n} \set{\pare{1+\alpha\norm{n}_q^2}^{-1} \psca{\tilde{f}_m^\perp,m}_q \tilde{f}_n^\perp - \pare{1+\alpha\norm{m}_q^2}^{-1} \psca{\tilde{f}_n^\perp,n}_q \tilde{f}_m^\perp}\\
	&\qquad + \frac{1}2 \sin\psca{m+n,x}_q P_{m+n} \set{\pare{1+\alpha\norm{n}_q^2}^{-1} \psca{\tilde{f}_m^\perp,m}_q \tilde{f}_n^\perp + \pare{1+\alpha\norm{m}_q^2}^{-1} \psca{\tilde{f}_n^\perp,n}_q \tilde{f}_m^\perp},
\end{align*}
and
\begin{align*}
	&\Bcal(g_m+g_n)\\ 
	&\quad = \frac{1}2 \sin\psca{m-n,x}_q P_{m-n} \set{\pare{1+\alpha\norm{n}_q^2}^{-1} \psca{\tilde{g}_m^\perp,m}_q \tilde{g}_n^\perp - \pare{1+\alpha\norm{m}_q^2}^{-1} \psca{\tilde{g}_n^\perp,n}_q \tilde{g}_m^\perp}\\
	&\qquad - \frac{1}2\sin\psca{m+n,x}_q P_{m+n} \set{\pare{1+\alpha\norm{n}_q^2}^{-1} \psca{\tilde{g}_m^\perp,m}_q \tilde{g}_n^\perp + \pare{1+\alpha\norm{m}_q^2}^{-1} \psca{\tilde{g}_n^\perp,n}_q \tilde{g}_m^\perp}.
\end{align*}
Then, we can repeat the above argument for $a$ to find $b$. Lemma \ref{le:Saturation} is proved. \qquad $\blacksquare$ 

\bigskip

In what follows, we recall that, for any finite-dimensional subspace $E$ of $V^3(\TT^2)^2$, we have defined $\mathcal{F}(E)$ as the largest vector subspace of $V^3(\TT^2)^2$ such that, for any $\overline{\eta} \in \mathcal{F}(E)$, there exist
$$k\in\NN^*;\; \alpha_1,\ldots, \alpha_k > 0;\; \eta, \rho^1, \ldots, \rho^k \in E$$ satisfying
$$\overline{\eta} = \eta - \sum_{j=1}^k \alpha_j \Bcal(\rho^j).$$
We have also defined a sequence of subspace $$E = E_0 \subset E_1 \subset \ldots \subset E_n \subset \ldots$$ such that, for any $n\in\NN$ we have $E_{n+1} = \Fcal(E_n)$ and we set $$E_\infty = \bigcup_{n=0}^\infty E_n.$$ The saturation property in Theorem \ref{th:FG2Saturation} can be made precise as follows
\begin{thm}[Saturation Property]
	\label{th:Saturation} Let $q = (q_1,q_2)$, $q_1, q_2 > 0$, let $E$ be a finite-dimensional subspace of $V^3(\TT^2)^2$ and for any $N\in\NN^*$, let $\Hcal^N_q$ be defined as in \eqref{eq:HN}. If $E \supset \Hcal^3_q$, then $E_\infty \supset \Hcal^N_q$, for any $N \in \NN$, $N \geq 3$. 
\end{thm}

\noindent\textbf{Proof. } Inspired by the argument of [\cite{S2006}, Theorem 2.5], we will prove Theorem \ref{th:Saturation} by recurrence that, if $E = E_0 \supset \Hcal^3_q$ then, for any $k \geq 0$, we have 
\begin{equation}
	\label{eq:E2k} E_{2k} \supset \Hcal^{k+3}_q.
\end{equation}
It is evident that for $k=0$, \eqref{eq:E2k} is true. Let $k \geq 1$. We suppose that, for any $k'\in\NN$, $0 \leq k' < k$, we have $E_{2k'} \supset \Hcal^{k'+3}_q$. In order to prove that \eqref{eq:E2k} is true for $k$, we only need to prove that, for any $l\in\ZZ^2 \setminus \set{0}$, $\abs{l} = k + 3$, we have $c_l, s_l \in E_{2k}$.

\vspace{0.3cm}

\noindent \textit{1. First case: } If $l = (l_1,l_2)$ with $l_1 \neq 0$ and $l_2 \neq 0$. In this case, since $\abs{l_1} + \abs{l_2} > 3$, without loss of generality, we can suppose that $l_1 \geq 2$. We choose $m = (l_1-1,l_2)$ and $n = (1,0)$. Then, we have
$$m + n = l, \quad \norm{m}_q > \norm{n}_q, \quad \abs{m} = k+2, \quad \abs{n} = 1, \quad \abs{m-n} = k+1$$ and $m$, $n$ are not parallel in $\RR^2$. Applying Lemma \ref{le:Saturation}, we obtain the existence of $$a, b \in Span\set{\Ccal_m,\Ccal_n,\Scal_m,\Scal_n} \subset \Hcal^{k+2}_q \subset E_{2k-2}$$ such that $$\Bcal(a) + c_l, \; \Bcal(b) + s_l \in Span\set{\Ccal_{m-n},\Scal_{m-n}}.$$ Thus, there exist $$f, g \in Span\set{\Ccal_{m-n},\Scal_{m-n}} \subset \Hcal^{k+1}_q \subset E_{2k-2}$$ such that $$c_l = f - \Bcal(a), \quad s_l = g - \Bcal(b).$$ Then, using the definition of $E_{2k-1} = \Fcal(E_{2k-2})$, we deduce that $$c_l, s_l \in E_{2k-1} \subset E_{2k}.$$

\vspace{0.3cm}

\noindent \textit{2. Second case: } If $l = (l_1,l_2)$ with $l_1 = \abs{l} = 3$ and $l_2 = 0$. In this case, we choose $m = (l_1-1,1)$ and $n = (1,-1)$. Then, we have
$$m + n = l, \quad \norm{m}_q > \norm{n}_q, \quad \abs{m} = k+3, \quad \abs{n} = 2, \quad \abs{m-n} = k+3$$ and $m$, $n$ are not parallel in $\RR^2$. Since all the components of the vectors $m$, $n$, $m-n$ are not zero, we can apply the result of the first case and we deduce that
$$Span\set{\Ccal_m,\Ccal_n,\Scal_m,\Scal_n} \subset E_{2k-1}$$ and $$Span\set{\Ccal_{m-n},\Scal_{m-n}} \subset E_{2k-1}.$$ Now, applying Lemma \ref{le:Saturation}, we obtain the existence of $$a, b \in Span\set{\Ccal_m,\Ccal_n,\Scal_m,\Scal_n} \subset E_{2k-1}$$ such that $$\Bcal(a) + c_l, \; \Bcal(b) + s_l \in Span\set{\Ccal_{m-n},\Scal_{m-n}} \subset E_{2k-1}.$$ As in the first case, we can deduce that
$$c_l, s_l \in \Fcal(E_{2k-1}) = E_{2k}.$$ Theorem \ref{th:Saturation} is proved. \qquad $\blacksquare$

%----- %----- %----- %----- %----- %

\bigskip
\section{Approximate controllability by high-mode controls reduced to approximate controllability by low-mode controls} \label{se:Dense}

The goal of this section is to prove the main theorem \ref{th:FG2approx} by proving that we can reduce the control of the system \eqref{sys:FG2eta} by high-mode controls to controls in $\Hcal^3_q$. Let $T > 0$, $\ee > 0$ and $u_0, u_T \in V^4(\TT^2)^2$. We set $$\Ucal_0 = u_0 - \alpha\DD u_0 \quad\mbox{and}\quad \Ucal_T = u_T - \alpha\DD u_T.$$ For any $t \in [0,T]$, let 
$$\overline{\Ucal}(t) = \frac{1}T \pare{I-\alpha\DD} \pare{(T-t)u_0 + tu_T}.$$ Then, $\overline{\Ucal}$ is solution of the system \eqref{sys:FG2etaCal} with
$$\overline{\Ucal}(0) = \Ucal_0 \quad\mbox{and}\quad \eta = \dd_t\overline{\Ucal} + \Lcal\overline{\Ucal} + \Bcal(\overline{\Ucal}) - \PP f.$$
It is clear that $$\overline{\Ucal} \in C(0,T,V^2(\TT^2)^2) \quad\mbox{and}\quad \eta \in L^{\infty}(0,T,V^1(\TT^2)^2).$$
Let $k \in \NN^*$, $k \geq 3$ and let $$\eta_k = P_k\pare{\dd_t\overline{\Ucal} + \Lcal\overline{\Ucal} + \Bcal(\overline{\Ucal}) - \PP f},$$ where $P_k$ is the projection onto $\Hcal^k_q$. Let $\overline{\Ucal}_k$ be the solution of the system \eqref{sys:FG2etaCal} with  
$$\overline{\Ucal}_k(0) = \Ucal_0 \quad\mbox{and}\quad \eta = \eta_k.$$
Applying Theorem \ref{th:FG2pertCal}, we can choose $k$ so large that 
\begin{equation*}
	\norm{\overline{\Ucal}_k(T) - \overline{\Ucal}(T)}_{V^1(\TT^2)^2} \leq C \norm{\eta_k - \pare{\dd_t\overline{\Ucal} + \Lcal\overline{\Ucal} + \Bcal(\overline{\Ucal}) - \PP f}}_{L^2(0,T,V^1(\TT^2)^2)} \leq \ee.
\end{equation*}
Now, Theorem \ref{th:Saturation} (see \eqref{eq:E2k}) implies that $\Hcal^k_q \subset E_N$, where $N = 2(k-3)$. Now, we set $\Ucal^N = \overline{\Ucal}_k$. Applying Theorems \ref{th:ConvexP} and \ref{th:ExtensionP}, we can contruct a sequence of controls $\eta_j \in E_j$, $j \in \set{0, \ldots, N}$, and a sequence of solutions $\Ucal^j$ of the system \eqref{sys:FG2etaCal} with  
$$\Ucal^j(0) = \Ucal_0 \quad\mbox{and}\quad \eta = \eta_j,$$ such that $$\norm{\Ucal^{j-1}(T) - \Ucal^j(T)}_{V^1(\TT^2)^2} \leq \frac{\ee}{2^{N-j}},$$ for any $j \in \set{1,\ldots,N}$. Thus, $\Ucal^0$ is the solution of the system \eqref{sys:FG2etaCal} with
$$\Ucal^0(0) = \Ucal_0 \quad\mbox{and}\quad \eta = \eta_0 \in E_0 = \Hcal^3_q,$$ and moreover, we have
\begin{equation*}
	\norm{\Ucal^0(T) - \Ucal_T}_{V^1(\TT^2)^2} \leq \ee + \sum_{j=1}^N \frac{\ee}{2^{N-j}} < 3\ee.
\end{equation*}
Finally, we set $u = (I-\alpha\DD)^{-1} \Ucal^0$. Then $u$ is the solution of the system \eqref{sys:FG2eta}, with 
$$u(0) = u_0 \quad\mbox{and}\quad \eta = \eta_0 \in \Hcal^3_q.$$ Moreover, we have
\begin{equation*}
	\norm{u(T) - u_T}_{V^3(\TT^2)^2} \leq \frac{1}{\min\set{1,\alpha}} \norm{\Ucal^0(T) - \Ucal_T}_{V^1(\TT^2)^2} \leq \frac{3\ee}{\min\set{1,\alpha}}.
\end{equation*}
Theorem \ref{th:FG2approx} is proved. \qquad $\blacksquare$

%\vspace{3cm}

\bigskip

\section*{Acknowledgments}

This work was partially done during the visit of the first author at the Vietnam Institute for Advanced Study in Mathematics (VIASM). The first author is thankful to the VIASM for the financial support and for the very kind hospitality of the institute and of all the staff.

\end{document}